\documentclass[11pt,letterpaper]{amsart}
\usepackage[margin=1in]{geometry}
\usepackage{graphicx,nccmath}
\usepackage{mathtools}
\usepackage{enumerate}
\usepackage{booktabs}
\usepackage{tabularx}
\usepackage{multirow}
\usepackage{soul}

\usepackage[scr=boondoxupr,cal=euler]{mathalpha}

\makeatletter
\newcommand{\distas}[1]{\mathbin{\overset{#1}{\kern\z@\sim}}}%
\newsavebox{\mybox}\newsavebox{\mysim}
\newcommand{\oversim}[1]{%
 \savebox{\mybox}{\hbox{\kern2pt$\scriptstyle#1$\kern2pt}}%
 \savebox{\mysim}{\hbox{$\sim$}}%
 \mathbin{\overset{#1}{\kern\z@\resizebox{\wd\mybox}{\ht\mysim}{$\sim$}}}%
}
\makeatother

\DeclareMathOperator{\sgn}{sgn}
\DeclareMathOperator{\var}{Var}
\DeclareMathOperator{\Nor}{N}
\DeclareMathOperator{\Cau}{Cauchy}
\DeclareMathOperator{\FN}{FoldedN}
\DeclareMathOperator{\Lev}{Levy}
\DeclareMathOperator{\Sta}{Stable}
\DeclareMathOperator{\Exp}{Exp}
\DeclareMathOperator{\Ber}{Bernoulli}
\DeclareMathOperator{\Bet}{Beta}
\DeclareMathOperator{\Gam}{Gamma}
\DeclareMathOperator{\Wb}{Weibull}
\DeclareMathOperator{\LN}{LogN}
\DeclareMathOperator{\Uni}{Unif}

\newcommand{\dto}{\xrightarrow{d}}
\newcommand{\deq}{\overset{d}{=}}

\usepackage{hyperref,xcolor}
\hypersetup{
colorlinks,
linkcolor={red!50!black},
citecolor={blue!50!black},
urlcolor={blue!80!black}
}

\newtheorem{theorem}{Theorem}[section]
\newtheorem{proposition}[theorem]{Proposition}
\newtheorem{lemma}[theorem]{Lemma}
\newtheorem{corollary}[theorem]{Corollary}
\theoremstyle{definition}
\newtheorem{definition}[theorem]{Definition}

\theoremstyle{remark}

\begin{document}

\title{Statistics of the Shepp $p$-product}
\author{Wenxuan Guo}
\author{Lek-Heng Lim}
\address{Department of Statistics, University of Chicago, Chicago, IL 60637}
\email{wxguo@uchicago.edu, lekheng@uchicago.edu}

\begin{abstract}
In 1962, Shepp famously discovered a product of normal random variables that preserves normality. The Shepp product, which takes the form $XY/\sqrt{X^2 + Y^2}$, has since been thoroughly studied and has found numerous connections to other areas of statistics. Among other things, it has an extension to $n$ normal variables, gives a multiplicative analogue of the central limit theorem, and applies rather unexpectedly to genomics as a test statistic for alignment-free sequence analysis. The Shepp product is evidently the $p =2$ special case of $XY/\sqrt[p]{X^p + Y^p}$ that we call the Shepp $p$-product. We will show that the Shepp $p$-product, particularly when $p = 1$ and $p = \infty$ (the latter in a limiting sense), is no less fascinating and applicable than the original $p = 2$ case. Just as the Shepp $2$-product preserves normal distributions, the Shepp $1$-product preserves Cauchy distributions, and the Shepp $\infty$-product preserves exponential distributions. In fact, the converse is also true in an appropriate sense, allowing us to characterize the Cauchy, normal, and exponential distributions as the unique distributions preserved by the Shepp $p$-product for $p =1,2,\infty$, respectively. We will study the multiplicative analogue of infinite divisibility with respect to the Shepp $p$-product, establish an asymptotic theory for the Shepp $p$-product of $n$ i.i.d.\ random variables, and estimate the rates of convergence in Kolmogorov distance. Some of our results are new even for the $p=2$ case. We will also discuss new applications of the Shepp $p$-product in statistics, computational biology, and statistical physics. In particular, our new test statistic for exponentiality based on the Shepp $\infty$-product is more powerful, both theoretically and empirically, than the rule-of-thumb Kolmogorov--Smirnov test based on the relative efficiency.

\end{abstract}


\maketitle

\section{Introduction}\label{sec:intro}

It is well-known that the straightforward product $XY$ of two independent normal variables $X$ and $Y$ is not normally distributed. On the other hand, if $X \sim \Nor(0,\sigma_1^2)$ and $Y\sim \Nor(0,\sigma_2^2)$ are independent, then a product of the form
\begin{equation}\label{eq:S2}
Z = \frac{XY}{\sqrt{X^2 + Y^2}}
\end{equation}
satisfies $Z \sim \Nor(0,\sigma^2)$ with
\[
\frac{1}{\sigma} = \frac{1}{\sigma_1} + \frac{1}{\sigma_2},
\]
a beautiful result first discovered by Shepp \cite{Shepp1962} and that has since found numerous different proofs \cite{Baringhaus1988, Quine1994, Reid1987}. Shepp's product has many fascinating properties: a converse holds in an appropriate sense \cite[Theorem~1]{Novak2007}; it may be extended to $n$ independent $X_i \sim \Nor(0,\sigma_i^2)$, $i =1,\dots,n$,
\begin{equation}\label{eq:genS}
S(X_1, \dots, X_n) \coloneqq \frac{X_1\cdots X_n}{\sqrt{\sum X_{i_1}^2\cdots X_{i_{n-1}}^2}}\sim \Nor\bigl(0, \sigma^2), \qquad \frac{1}{\sigma} =\frac{1}{\sigma_1} + \dots + \frac{1}{\sigma_n},
\end{equation}
where the sum in the denominator is over all $1 \le i_1 <i_2 <\dots<i_{n-1} \le n$; and in addition it satisfies a central limit theorem: For \emph{any} i.i.d.\ random variables $X_i$, $i \in \mathbb{N}$, i.e., not necessarily normal, $n S(X_1, \dots, X_n)$ is asymptotically normally distributed as $n \to \infty$ \cite{Quine1999}.

Shepp's product also has various applications, and not just in statistics, but in areas like computational biology and statistical physics too. We briefly review some of these as background for our proposed applications in this article. In statistics, the Shepp product $S$ may be used to construct scale-invariant normality tests \cite{Volkova2009}: For i.i.d.\ $X_1,\dots,X_n$, let $U_n$ and $V_n$ be the $U$- and $V$-statistical empirical distribution functions associated with $S$, i.e.,
\[
U_{n}(t) =\frac{2}{n(n-1)} \sum_{1 \le i<j \le n} \mathbb{I}_{\{S(X_i,X_j)<t/2\}},\qquad V_{n}(t)=\frac{1}{n^2} \sum_{i, j=1}^{n} \mathbb{I}_{\{S(X_i, X_j)<t/2\}},
\]
and let $F_n$ denote the empirical distribution function. Then the \emph{integral} and \emph{Kolmogorov statistics} may be defined by
\[
\mathit{IU}_n \coloneqq \int_{\mathbb{R}} \bigl(V_n(t) - F_n(t)\bigr) \, dF_n(t), \qquad
\mathit{DU}_n \coloneqq \sup_{t \in \mathbb{R}}\; \lvert U_{n}(t)-F_{n}(t)\rvert
\]
respectively. In particular, $\mathit{IU}_n$ is \emph{scale-free}, unlike the Kolmogorov--Smirnov test. The asymptotic behaviors and the power of these test statistics may be found in \cite{Volkova2009}.

In genomics, a statistic based on Shepp product is used to test the difference between two cDNA sequences in alignment-free sequence analysis \cite{Reinert2009}:
\[
D_2^S =\sum_{w \in \mathcal{A}^k } \frac{X_w Y_w}{\sqrt{\smash[b]{X_w^2 +Y_w^2}}}.
\]
Here, $X_w$ and $Y_w$ are centralized count variables measuring the frequency of word $w$ in the two sequences, and $\mathcal{A}^k$ denotes the set of words of length $k$ in the alphabet $\mathcal{A}$. Under some mild conditions, $D_2^S$ converges to the normal distribution under the null hypothesis. This has been extended to the $n$-variate Shepp product \eqref{eq:genS} for multiple sequence comparison and is used to conduct extensive simulations \cite{Ren2013}. We refer interested readers to \cite{Song2013} for a review of the role of Shepp's product in alignment-free sequence analysis.

The Shepp product is also a special case of a \emph{Wald statistic} under a singular hypothesis \cite{Glonek1993, Drton2016}. Suppose we have estimators $\hat{\theta}_1$ and $\hat{\theta}_2$ for unknown parameters $\theta_1$ and $\theta_2$ respectively. We are interested in testing $H_0: \theta_1\theta_2 = 0$ under the singular condition $\theta_1 = \theta_2 = 0$. If we have $\sqrt{n}\bigl((\hat{\theta}_1, \hat{\theta}_2) - ({\theta}_1, {\theta}_2)\bigr) \dto \Nor(0, I_2)$, a Wald statistic can be constructed to test $H_0$. Surprisingly, it converges to the Shepp product given $\theta_1 = \theta_2 = 0$:
\[
  \text{Wald statistic} \coloneqq \sqrt{n}\frac{\hat{\theta}_1\hat{\theta}_2}{\sqrt{\smash[b]{\hat{\theta}_1^2 + \hat{\theta}_2^2}}} \dto \frac{XY}{\sqrt{X^2+Y^2}},\qquad X, Y\oversim{iid}\Nor(0, 1).
\]
Because of the singular hypothesis, the asymptotic distribution of the Wald statistic is no longer $\Nor(0, 1)$ but $\Nor(0, \frac{1}{4})$, and, indeed, this follows from the distribution theory of the Shepp product. A consequence is that the standard Wald test is conservative when the true parameters are zero.

The goal of our article is to show that the Shepp $p$-product, defined by
\begin{equation}\label{eq:Sp}
S_p(X, Y) \coloneqq
\begin{cases}
\dfrac{X Y}{(X^p + Y^p)^{1/p}} & 1 \le p < \infty,\\[3ex]
\dfrac{1}{\max(1/X, 1/Y)} & p = \infty,
\end{cases}
\end{equation}
is no less fascinating and useful than Shepp's original product in \eqref{eq:S2}, which corresponds to the case $p =2$ and will henceforth be called the Shepp $2$-product. We should perhaps mention that while it is the point-of-view undertaken in this article, it is unusual for the expression in \eqref{eq:S2} to be regarded as a product of random variables. Nevertheless, this perspective is convenient as we may easily show that the Shepp $p$-product in \eqref{eq:Sp} is associative for all $p \in [1,\infty]$ (it is obviously commutative); this in turn allows us to unambiguously define a Shepp $p$-product of $n$ random variables $S_p(X_1, \dots, X_n)$, which agrees with \eqref{eq:genS} when $p=2$. We will study the Shepp $p$-product with an emphasis on its asymptotic behavior and discuss several applications, particularly for the cases $p=1$ and $\infty$. Observe that for $p \in [1,\infty)$,
\begin{equation}\label{eq:transform}
\frac{1}{S_p(X_1, \dots, X_n)^p} =  \sum_{i = 1}^n \frac{1}{X_i^p} \quad\text{and}\quad
\frac{1}{S_\infty(X_1, \dots, X_n)} = \max_{i = 1,\dots, n}\Bigl(\frac{1}{X_i}\Bigr),
\end{equation}
which connects Shepp $p$-product to sum of reciprocals for finite $p$, and extreme order statistics for $p = \infty$. The latter may be analyzed using standard distribution theory, but the former is not well-studied --- reciprocals are typically heavy-tailed, and their i.i.d.\ sum does not exhibit the usual asymptotic normality with $n^{-\frac12}$-scaling. Our Shepp $p$-product opens a door for the study of such sums using the theory of stable distribution for finite $p$, and the theory of extreme order statistics for $p = \infty$. We summarize some of our key findings and several new applications below.

\subsection*{Characterization of Cauchy, normal, and exponential distributions}

In Section~\ref{sec:infdiv}, we show that the Shepp $1$-product is to the  Cauchy distribution and the Shepp $\infty$-products is to the exponential distributions what the Shepp $2$-product is to the normal distribution. Let $\gamma,\lambda \in \mathbb{R}$ with $\lambda > 0$. We write $\Cau(\gamma, \lambda)$ for the Cauchy distribution with location and scale parameters $\gamma$ and $\lambda $, respectively, i.e., whose distribution function is
$F(x) = \pi^{-1} \arctan \bigl((x-\gamma)/\lambda\bigr)+1/2$. We write $\Exp(\lambda)$ for the exponential distribution with rate parameter $\lambda$, i.e., whose distribution function is $F(x) = 1 - \exp(-\lambda x)$.

If $X \sim \Cau(\gamma_1, \lambda_1)$ and $Y \sim \Cau(\gamma_2, \lambda_2)$ are independent, then their Shepp $1$-product
\begin{gather}\label{eq:S1}
Z=\frac{X Y}{X + Y} \sim  \Cau(\gamma, \lambda),
\shortintertext{where}
\begin{medsize}
\gamma = \dfrac{\dfrac{\gamma_1}{\gamma_1^2 + \lambda_1^2}+ \dfrac{\gamma_2}{\gamma_2^2 + \lambda_2^2}}{\biggl[\dfrac{\gamma_1}{\gamma_1^2 + \lambda_1^2} + \dfrac{\gamma_2}{\gamma_2^2 + \lambda_2^2}\biggr]^2 + \biggl[\dfrac{\lambda_1}{\gamma_1^2 + \lambda_1^2} + \dfrac{\lambda_2}{\gamma_2^2 + \lambda_2^2}\biggr]^2}, \quad
 \lambda = \dfrac{\dfrac{\lambda_1}{\gamma_1^2 + \lambda_1^2}+ \dfrac{\lambda_2}{\gamma_2^2 + \lambda_2^2}}{\biggl[\dfrac{\gamma_1}{\gamma_1^2 + \lambda_1^2}+\dfrac{\gamma_2}{\gamma_2^2 + \lambda_2^2}\biggr]^2 + \biggl[\dfrac{\lambda_1}{\gamma_1^2 + \lambda_1^2}+\dfrac{\lambda_2}{\gamma_2^2 + \lambda_2^2}\biggr]^2}.
\end{medsize}\nonumber
\end{gather}
Likewise, if $X \sim \Exp(\lambda_1)$ and $Y \sim \Exp(\lambda_2)$ are independent, then their Shepp $\infty$-product
\begin{equation}\label{eq:Sinfty}
Z=\frac{1}{\max(1/X, 1/Y)} \sim  \Exp(\lambda_1 + \lambda_2).
\end{equation}
These are clearly analogues of \eqref{eq:S2}. More importantly, the converse of the above statements also holds in an appropriate sense, i.e., the products in \eqref{eq:S1} and \eqref{eq:Sinfty} characterize the Cauchy and exponential distributions in the following manner. Let $p \in [1,\infty]$.
We say that the Shepp $p$-product uniquely characterizes a distribution if it is  the unique distribution where
\begin{equation}\label{eq:distri}
c S_p(X, Y) \deq X, \quad c > 0
\end{equation}
holds for any i.i.d.\ $X_1$, $X_2$ in that distribution, and for some constant $c$. Here and subsequently, we write $\deq$ to mean equality in distribution.

Under a mild assumption that the distribution function has \emph{regular variation index} $r=1$ (formally defined in Section~\ref{sec:char} and \cite[Chapter~1]{Seneta1976} but roughly it means that the density is continuous and nonzero at zero), the Shepp $2$-product is known to uniquely characterize the normal distribution with $c = 2$ \cite{Bansal1999, Galambos2003, Galambos2004, Hamedani2007, Novak2007}. This is in fact the basis for the normality test based on Shepp $2$-product \cite{Volkova2009} described earlier. We will show in Section~\ref{sec:char} that under the same assumption, the Shepp $1$-product uniquely characterizes the Cauchy distribution and the Shepp $\infty$-product uniquely characterizes the exponential distribution. These results are in turn used in Section~\ref{sec:test} to construct scale-invariant tests for Cauchy and exponential distributions.
We will also go beyond the $r=1$ case and completely determine the distributions uniquely characterized by the Shepp $1$- and $\infty$-product for all possible values of $r$, as well as those uniquely characterized by the Shepp $2$-product for $0 < r < 2$.

\subsection*{Infinite divisibility}

The notion of infinite divisibility is an important tool in the study of central limit theorems and characteristic functions \cite{Gnedenko1954, Lukacs1960}. It has a multiplicative analogue for the Shepp $p$-product, namely, a probability distribution is \emph{Shepp $p$-infinitely divisible} if it can be expressed as a distribution of a Shepp $p$-product of an arbitrary number of i.i.d.\ random variables. In Section~\ref{sec:infdiv}, we will study a wide range of Shepp $p$-infinitely divisible distributions for $p = 1, \infty$ as well as any even-integer-valued $p \in 2\mathbb{N}$. For finite values of $p$, our result essentially follows from the infinite divisibility of stable distributions. For $p = \infty$, we will show that a family of exponential-like distributions (to be formally defined in Section~\ref{sec:infinfdiv}, but roughly they are some power of exponential distributions), including the exponential distribution itself, are Shepp $\infty$-infinitely divisible distributions.

\subsection*{Asymptotic theory}

We will establish a central limit theorem for the Shepp $p$-product of any sequence of i.i.d.\ random variables $X_1,X_2,X_3,\dots$ for $p =1,\infty$ or $p \in 2\mathbb{N}$ in Section~\ref{sec:clt}. We will show that the asymptotic distribution of $S_p(X_1,\dots, X_n)$ is Cauchy if $p = 1$ and is a variant of the exponential distribution if $p = \infty$. Our result for even-integer-valued $p$ generalizes the asymptotic normality result for $p=2$ in \cite{Galambos2003, Quine1999}.

\subsection*{Rates of convergence}

In the second half of Section~\ref{sec:clt}, we quantify the convergence rates in the Kolmogorov distance, and thereby derive Berry--Esseen-type bounds for the Shepp $p$-product  for $p =1,\infty$ and $p \in 2\mathbb{N}$. A departure from classical central limit theorem and Berry--Esseen bound is that for the Shepp $p$-product, the value of the regular variation index $r$ has a non-negligible effect on the convergence. We will apply the notion of \emph{domain of normal attraction}  \cite[Section~35]{Gnedenko1954} in stable distribution theory to extreme value distributions; these are, to the best of our knowledge, the first results of the convergence rate of the domain of normal attraction of extreme value distributions.

\subsection*{Exponentiality and Cauchy tests}

In Section~\ref{sec:test}, we use our distributional characterizations for the Shepp $\infty$- and $1$-products in Section~\ref{sec:char} to construct exponentiality and Cauchy tests, defining a scale-invariant integral statistic $\mathit{IU}_n$ and a Kolmogorov-type statistic $\mathit{DU}_n$ analogous to \cite{Volkova2009}. We determine the asymptotic distributions under the null hypotheses to ascertain the validity of our testing procedures. While one of our statistics has an explicit asymptotic distribution, the others do not and for these, we discuss simulation methods for computing the critical values for the corresponding tests. With regards to the power, we adopt the notion of local Bahadur efficiency to analyze exponentiality tests. Our tests are proved to be more efficient than the Kolmogorov--Smirnov test, a conclusion that is further confirmed through extensive simulations.

\subsection*{Alignment-free sequence analysis}

Following \cite{Reinert2009} we propose several new test statistics in Section~\ref{sec:align} for alignment-free sequence analysis inspired by the Shepp $p$-product:
\[
D_p^S \coloneqq \sum_{w \in \mathcal{A}^k } S_p(X_w, Y_w), \qquad
C_1^S \coloneqq \sum_{w, w^\prime \in \mathcal{A}^k }
\frac{X_w X_{w^\prime}}{X_wY_{w^\prime} + X_{w^\prime}Y_w},
\]
where again $X_w$'s and $Y_w$'s are centralized count variables and $\mathcal{A}^k$ denotes the set of words of length $k$ in the alphabet $\mathcal{A}$. Note that $C_1^S \ne D_1^S$.  We will determine the asymptotic distributions under the null hypotheses as the sequence length goes to infinity.

\subsection*{Conductivity of random resistor network}

Central limit theorems with hierarchical structure play an important role in statistical physics, particularly in the study of conductance in random media \cite[Section~4.2]{Chen2011}. We will discuss in Section~$\ref{sec:media}$ an application of the Shepp $1$-product to a random resistor network with resistors connected in series, providing various asymptotic distributions for the conductance of the network by drawing on our results in Section~\ref{sec:clt}.

\section{Basic properties of the Shepp $p$-product}

We will begin by stating the Shepp $p$-product formally.
\begin{definition}[Shepp $p$-product]
For $p\in [1, \infty]$, the \emph{Shepp $p$-product} of two random variables $X$ and $Y$ is given by $S_p(X, Y)$ as defined in \eqref{eq:Sp}.
\end{definition}
The value $S_p(X, Y)$ is not well-defined for $X$ and $Y$ having point mass at zero; to avoid this we shall always assume that $\mathbb{P}(X = 0) = \mathbb{P}(Y = 0) = 0$. In addition, to avoid having to discuss a complex-valued Shepp $p$-product or restricting ourselves to positive-valued random variables, we will henceforth assume that $p$ is either $1$, $\infty$, or an even integer. Nevertheless, we do not wish to discount the possibility that the Shepp $p$-product as defined in \eqref{eq:Sp} may have interesting properties and applications for other values of $p \in 2 \mathbb{N} + 1$, just that such explorations will have to be deferred to future work.

The Shepp $p$-product is commutative and associative. Consequently, we may unambiguously define the Shepp $p$-product of any $n$ random variables by recursion and will not need to pay attention to the order of factors in a Shepp $p$-product.

\begin{proposition}[Commutativity and Associativity]\label{prop:assoc}
Let $X, Y, Z$ be real-valued random variables and $p \in \{1,\infty\} \cup 2\mathbb{N}$. Then
\[
S_p(X,Y) = S_p(Y,X), \qquad S_p\bigl(S_p(X, Y), Z\bigr) = S_p\bigl(X, S_p(Y, Z)\bigr).
\]
\end{proposition}
\begin{proof}
Commutativity is immediate.
We write $\sgn(0)=0$ and $\sgn(x) =x/\lvert x \rvert$ for $ x\in \mathbb{R}\setminus \{0\}$.  For $p<\infty$,
\[
S_p(X, Y)= \frac{\sgn(XY)}{\bigl(\sgn(Y)^p/|X|^p + \sgn(X)^p/|Y|^p\bigr)^{1/p}}.
\]
Therefore $S_1(X, Y)= 1/(1/X + 1/Y)$ and $S_p(X, Y)= \sgn(XY)/(1/|X|^p + 1/|Y|^p)^{1/p}$ for $p\in 2\mathbb{N}$. From
\[
 S_1\bigl(S_1(X, Y), Z\bigr) = \frac{1}{1/X + 1/Y + 1/Z}, \quad
 S_p\bigl(S_p(X, Y), Z\bigr) = \frac{\sgn(XYZ)}{(1/|X|^p + 1/|Y|^p + 1/|Z|^p)^{1/p}},
\]
we see that associativity holds for $p = 1$ and $p \in 2\mathbb{N}$. For $p = \infty$, this follows from
\[
 S_\infty\bigl(S_\infty(X, Y), Z\bigr) = \frac{1}{\max(1/X, 1/Y, 1/Z)}. \qedhere
\]
\end{proof}

With associativity, we may define an $n$-fold Shepp $p$-product by
\begin{equation}\label{eq:recur}
S_p(X_1, \dots, X_n) = S_p\bigl(S_p(X_1 , \dots, X_{n-1}), X_n\bigr)
\end{equation}
or indeed in any arbitrary order. The proof of Proposition~\ref{prop:assoc} also shows that for $p < \infty$,
\[
 S_p(X_1, \dots, X_n) = \frac{X_1\cdots X_n}{\bigl(\sum X_{i_1}^p\cdots X_{i_{n-1}}^p\bigr)^{1/p}}
\]
where the sum is over all $1 \le i_1 <i_2 <\dots<i_{n-1} \le n$, and for $p = \infty$,
\[
 S_\infty(X_1, \dots,X_n) = \frac{1}{\max(1/X_1, \dots, 1/X_n)}.
\]

The Shepp $p$-product as a stochastic process indexed by $n$, the number of variables, is a Markov process.
\begin{theorem}[Markov property]
Let $(\Omega, \mathcal{F}, \mathbb{P})$ be a probability space and $X_i:\Omega\to\mathbb{R}$ be independent random variables, $i\in\mathbb{N}$. We define the \emph{Shepp $p$-process} by
\[
Y_n \coloneqq S_p(X_1, \dots, X_n),
\]
adapted to the filtration $\mathcal{F}_n\coloneqq \sigma(X_1, \dots, X_n)$ for each $n \in \mathbb{N}$. Then $\{Y_n :  n \in \mathbb{N}\}$ is a Markov process.
\end{theorem}

\begin{proof}
For every nonnegative Borel measurable $f:\mathbb{R}\to\mathbb{R}$ and for every $m\le n$, we have $\mathbb{E}\bigl(f(Y_n)\bigm| \mathcal{F}_m\bigr) = \mathbb{E}\bigl(f\bigl(S_p(Y_m, X_{m+1}, \dots, X_n)\bigr)\bigm| \mathcal{F}_m\bigr)$ 
by \eqref{eq:recur}.
Since $S_p :\mathbb{R}^{n - m + 1}\to\mathbb{R}$ is Borel measurable, $\Tilde{f}\coloneqq  f*S_p$ is a positive Borel measurable function. We have $\mathbb{E}\bigl(f(Y_n)\bigm| \mathcal{F}_m\bigr) = \mathbb{E}\bigl(\Tilde{f}(Y_m, X_{m+1}, \dots, X_n)\bigm| \mathcal{F}_m\bigr)$. 
Since $Y_m$ is $\mathcal{F}_m$-measurable and $X_{m+1}, \dots, X_{n}$ are independent of $\mathcal{F}_m$, we may deduce from \cite[Lemma~A.3]{Schilling2012} that $\mathbb{E}\bigl(\Tilde{f}(Y_m, X_{m+1}, \dots, X_n)\bigm| \mathcal{F}_m\bigr) = g(Y_m)$, 
where $g(y) = \mathbb{E}\bigl(\Tilde{f}(y, X_{m+1},\dots, X_{n})\bigr)$. Hence, $\{Y_n :  n \in \mathbb{N}\}$ is a Markov process.
\end{proof}

To illustrate the benefits of the above perspective, we end this section with a note about the multivariate Shepp $2$-product and Wald statistic that extends results in \cite{Drton2016,Glonek1993}. In the following and throughout this article, we write $\dto$ to mean weak convergence, i.e., convergence in distribution.  Given $k$ parameters $\theta_1, \dots, \theta_k$ and a null hypothesis
\[
H_0:  \theta_1\theta_2\cdots\theta_k = 0,
\]
we let $(\widehat{\theta}_1,\dots, \widehat{\theta}_k)$ be an asymptotically normal estimator, i.e.,
\begin{equation}\label{eq:consistent}
 \sqrt{n}\bigl((\widehat{\theta}_1, \dots, \widehat{\theta}_k) - (\theta_1,\dots, \theta_k)\bigr)\dto \Nor(0, I_k).
\end{equation}
Then the corresponding Wald statistic is
\[
W_n=\sqrt{n}\frac{\widehat{\theta}_1 \cdots \widehat{\theta}_k}
 {\sqrt{\sum \widehat{\theta}^2_{i_1}\cdots\widehat{\theta}^2_{i_{k - 1}}}} = \sqrt{n} S_2(\widehat{\theta}_1, \dots, \widehat{\theta}_k)
\]
where the sum is over all $1 \le i_1 <i_2 <\dots<i_{k-1} \le k$; its asymptotic normality follows from standard theory. However, when $\theta_1 = \dots = \theta_k = 0$, it becomes a singular hypothesis \cite{Drton2016} and the standard asymptotic result no longer holds. Fortunately, we may  still derive the asymptotic distribution of $W_n$ based on \eqref{eq:consistent} and properties of the Shepp $2$-product,

\[
 W_n \dto \frac{X_1\cdots X_k}{\sqrt{\sum X_{i_1}^2\cdots X_{i_{k-1}}^2}}=S_2(X_1,\dots,X_k), \qquad (X_1, \dots, X_k)\sim\Nor(0, I_k),
\]
where the sum is again over all $1 \le i_1 <i_2 <\dots<i_{k-1} \le k$.
Since the asymptotic distribution has the form of a $k$-variate Shepp $2$-product, by either \cite{Quine1994} or results in Section~\ref{sec:2infdiv}, we deduce that $W_n\dto \Nor(0, 1/k^2)$. For $k=2$, we obtain the results in \cite[Section~2]{Drton2016} and \cite{Glonek1993}.

\section{Shepp $p$-infinite divisibility}\label{sec:infdiv}

We introduce a multiplicative analogue of infinite divisibility with respect to the Shepp $p$-product and show that some well-known distributions are Shepp $p$-infinitely divisible for various values of $p$.
\begin{definition}[Shepp $p$-infinite divisibility]\label{def:infdiv}
A probability distribution $\mathscr{P}$ is said to be \emph{Shepp $p$-infinitely divisible} if, for every positive integer $n$, there exist i.i.d.\ random variables $X_1, \dots, X_n$ whose Shepp $p$-product $S_p(X_1, \dots, X_n)$ has the same distribution $\mathscr{P}$. We say that a random variable is Shepp $p$-infinitely divisible if its distribution has the property.
\end{definition}
We begin by showing that Definition~\ref{def:infdiv} is not vacuous by exhibiting a Shepp $p$-infinitely divisible distribution for any $p \in (0, \infty)$ in Section~\ref{sec:pinfdiv}. For the special cases $p = 1,\infty$, $p \in 2\mathbb{N}$, we will give, in Sections~\ref{sec:1infdiv}, \ref{sec:2infdiv}, \ref{sec:infinfdiv}, respectively, the exact relations between the parameters of the distributions of $X_1,\dots,X_n$ and those of their Shepp $p$-product, thereby deducing Shepp $p$-infinite divisibility of these distributions.

For easy reference, we briefly remind readers of the notion of stable distributions before discussing results about Shepp $p$-infinite divisibility.

\subsection{Stable distributions}

A random variable $X$ is said to have a \emph{stable distribution} if for any $a > 0$ and $b > 0$, there is a $c > 0$ and a $d \in \mathbb{R}$ such that
\[
aX_1 + bX_2 \deq cX+ d, 
\]
where $X_1,X_2$ are independent copies of $X$. Stable distributions play an important role in central limit theorems, the normal distribution being a standard example. We refer readers to \cite[Chapter~1]{Samorodnitsky1994} for more information. The standard definition above is less useful for us than the following characterization adapted from \cite[Section~34]{Gnedenko1954}.
\begin{theorem}[Gnedenko--Kolmogorov]\label{thm:GK}
A distribution is stable if and only if its characteristic function $\varphi$ can be represented by
\begin{equation}\label{eq:GK}
\log \varphi(t) =
\begin{cases}
i \gamma t-\lambda|t|^{\alpha}\Bigl[1+i \beta \tan \Bigl(\dfrac{\pi \alpha}{2}\Bigr)\sgn(t) \Bigr] & \alpha\in (0, 1)\cup(1, 2],\\[3ex]
i \gamma t-\lambda|t|\Bigl[1+i\dfrac{2\beta}{\pi} \sgn(t) \log |t|\Bigr] & \alpha = 1,
\end{cases}
\end{equation}
with $\gamma\in\mathbb{R}$, $\lambda\in [0, \infty)$, $\beta\in[-1, 1]$, and $\alpha\in (0, 2]$. We will henceforth refer to \eqref{eq:GK} as the canonical representation and denote the corresponding stable distribution by $\Sta(\alpha, \beta, \lambda, \gamma)$.
\end{theorem}
In other words, throughout our article, a stable distribution $\Sta(\alpha, \beta, \lambda, \gamma)$ is taken to be one whose characteristic function is given by \eqref{eq:GK}. It includes as notable special cases (see \cite[Section~34]{Gnedenko1954}):
\begin{enumerate}[\upshape (i)]
\item\label{it:L} L\'evy distribution $\Lev(\gamma, \lambda)$ when $\alpha = 1/2$, $\beta = -1$;
\item\label{it:C} Cauchy distribution  $\Cau(\gamma, \lambda)$ when $\alpha = 1$, $\beta = 0$;
\item\label{it:N} normal distribution $\Nor(\gamma,2\lambda)$ when $\alpha = 2$.
\end{enumerate}
In \eqref{it:L} and \eqref{it:C}, $\gamma$ and $\lambda$ are often called the location and scale parameters, respectively. Note that $\beta$ is absent from \eqref{it:N} as the term involving $\beta$ in \eqref{eq:GK} vanishes when $\alpha =2$.

\subsection{Shepp $p$-infinitely divisible distributions}\label{sec:pinfdiv}

For any $p \in [1, \infty)$, a Shepp $p$-infinitely divisible distribution can be obtained as a power transform of a folded normal distribution. Specifically, if $X\sim\Nor(0,\sigma^2)$, then $Y=|X|\sim\FN(0,\sigma^2)$, and the distribution of $Y^{2/p}$ is Shepp $p$-infinitely divisible.

\begin{proposition}[Existence of Shepp $p$-infinitely divisible distributions]\label{prop:exist}
Let $p\in [1, \infty)$ and $n \in \mathbb{N}$.
Let $X_1, \dots, X_n$ be independent positive-valued random variables with
\[
X_i^{p/2} \sim \FN(0, \sigma_i^2).
\]
Let $T_p\coloneqq S_p(X_1, \dots, X_n)$. Then
\[
T_p^{p/2}\sim \FN(0, \sigma^2)
\]
with $1 / \sigma  = \sum_{i = 1}^n 1 / \sigma_i$.

\end{proposition}

\begin{proof}
By the definition of the Shepp $p$-product \eqref{eq:Sp},
\[
T_p = \frac{1}{(1/X_1^p + \dots + 1/X_n^p)^{1/p}},
\]
and thus $1/T_p^p = 1/X_1^p + \dots +  1/X_n^p$ with
$X_i^{p/2} \sim \FN(0, \sigma_i^2)$. Since $1/X_i^p$ follows the inverse square of a normal distribution, it follows the L\'evy distribution with the log characteristic function
\[
\log \varphi_i(t) = - \frac{1}{\sigma_i}|t|^{1/2}\bigl(1 - i\sgn(t)\bigr).
\]
Now $1/T_p^p$ has the log characteristic function
\[
\log \varphi(t) = - \sum_{i = 1}^n\frac{1}{\sigma_i}|t|^{1/2}\bigl(1 - i\sgn(t)\bigr).
\]
As $T_p$ and $1/T_p^p$ are positive random variables, we infer from the distribution of $1/T_p^p$ that
\[
T_p^{p/2}\sim \FN(0, \sigma^2)
\]
with $1/\sigma = 1/\sigma_1 + \dots + 1/\sigma_n$.
\end{proof}
Consequently, the distribution of $Y^{2/p}$ where $Y\sim\FN(0,\sigma^2)$ is Shepp $p$-infinitely divisible, for any $p \in [1,\infty)$ and $\sigma>0$. Indeed, $S_p(X_1,\dots,X_n)^{p/2} \sim\FN(0,\sigma^2)$ for any independent $X_1^{p/2}, \dots,  X_n^{p/2} \sim \FN(0,n^2\sigma^2)$.

\subsection{Shepp $1$-infinite divisibility}\label{sec:1infdiv}

We next exhibit a Shepp $1$-infinitely divisible distribution. A stable distribution $\Sta(\alpha, \beta, \lambda, \gamma)$ is as defined in \eqref{eq:GK}.

\begin{proposition}\label{prop:1param}
Let $\alpha\in(0, 2]$, $n \in \mathbb{N}$, and $X_1,\dots,X_n$ be independent with
$1/X_i \sim \Sta(\alpha, \beta_i, \lambda_i, \gamma_i)$, $i=1,\dots,n$. Let
$T_1 \coloneqq S_1(X_1, \dots, X_n)$. Then
\[
 \frac{1}{T_1} \sim \Sta(\alpha, \beta, \lambda, \gamma)
\]
with parameters $\beta, \lambda, \gamma$ given by
\begin{equation}\label{eq:1param}
 \lambda = \lambda_1 + \dots + \lambda_n,\quad
 \gamma = \gamma_1 + \dots + \gamma_n,\quad
 \beta = \frac{\lambda_1\beta_1 + \dots + \lambda_n\beta_n}{\lambda_1 + \dots + \lambda_n}.
\end{equation}
\end{proposition}
\begin{proof}
As stable distributions do not have a point mass at zero \cite[Chapter~1]{Samorodnitsky1994}, $1/X_i$ is well-defined as a random variable, and we have $1/X_1 + \dots + 1/X_n \sim \Sta(\alpha, \beta, \lambda, \gamma)$, with parameters as in \eqref{eq:1param}.
\end{proof}
Restricting Proposition~\ref{prop:1param} to i.i.d.\ $X_i$'s shows that the reciprocal of any stable random variable is Shepp $1$-infinitely divisible.

\begin{corollary}\label{cor:1infdiv}
The random variable $T_1$ such that $1 / T_1 \sim \Sta(\alpha, \beta, \lambda, \gamma)$ is Shepp $1$-infinitely divisible, i.e., for any $n \in \mathbb{N}$,
we have
\[
T_1\deq S_1(X_1, \dots, X_n)
\]
with independent 
 $1 / X_1, \dots, 1 / X_n \sim \Sta(\alpha, \beta, \lambda/n, \gamma/n)$.
\end{corollary}

Two special cases of Corollary~ \ref{cor:1infdiv} deserve highlight: the Cauchy distribution and the chi-square distribution.
\begin{corollary}[Cauchy distribution is Shepp $1$-infinitely divisible]\label{cor:Cparam}
Let $X_i \sim \Cau(\gamma_i, \lambda_i)$, $i =1, \dots, n$, be independent. Then $S_1(X_1, \dots, X_n)\sim \Cau(\gamma, \lambda)$ with
\[
\gamma = \dfrac{\displaystyle\sum_{i = 1}^n \dfrac{\gamma_i}{\gamma_i^2 + \lambda_i^2}}{\biggl[ \displaystyle\sum_{i = 1}^n \dfrac{\gamma_i}{\gamma_i^2 + \lambda_i^2}\biggr]^2 + \biggl[\displaystyle\sum_{i = 1}^n \dfrac{\lambda_i}{\gamma_i^2 + \lambda_i^2} \biggr]^2}, \qquad
 \lambda = \dfrac{\displaystyle\sum_{i = 1}^n \dfrac{\lambda_i}{\gamma_i^2 + \lambda_i^2}}{\biggl[ \displaystyle\sum_{i = 1}^n \dfrac{\gamma_i}{\gamma_i^2 + \lambda_i^2}\biggr]^2 + \biggl[\displaystyle\sum_{i = 1}^n \dfrac{\lambda_i}{\gamma_i^2 + \lambda_i^2} \biggr]^2}.
\]
In particular, $\Cau(\gamma, \lambda)$ is Shepp $1$-infinitely divisible with
\[
S_1(X_1, \dots, X_n) \sim \Cau(\gamma, \lambda)
\]
whenever $X_1,\dots,X_n \sim \Cau(n\gamma, n\lambda)$ are independent.
\end{corollary}
\begin{proof}
Observe that if $1 / X \sim \Sta(1, 0, \lambda, \gamma)$, then $X \sim \Cau\bigl(\gamma / (\gamma^2+\lambda^2), \lambda / (\gamma^2+\lambda^2)\bigr)$.
\end{proof}

In the following, we denote the chi-square distribution with $k$ degrees of freedom by $\chi^2_k$. To avoid clutter, we use $\chi^2_k$ instead of $\chi^2_k(0)$, dropping the noncentrality parameter.
\begin{corollary}[Chi-square distribution is Shepp $1$-infinitely divisible]\label{cor:chisq}
Let $\lambda_i^2 X_i \sim \chi_1^2$ for constants $\lambda_i > 0$, $i=1,\dots,n$, and let $X_1, \dots, X_n$ be independent. Then $\lambda^2 S_1(X_1, \dots, X_n)\sim \chi_1^2$ with $ \lambda \coloneqq \lambda_1 + \dots + \lambda_n$. 
In particular, $\chi_1^2/\lambda^2$ is Shepp $1$-infinitely divisible with
\[
\lambda^2 S_1(X_1, \dots, X_n) \sim \chi_1^2
\]
whenever $\lambda^2 X_1/n^2,\dots, \lambda^2 X_n/n^2 \sim \chi_1^2$ are independent.
\end{corollary}
\begin{proof}
Observe that if $1 / X \sim \Sta(1/2, -1, \lambda, 0 )$, then $\lambda^2 X \sim \chi^2_1$.
\end{proof}

\subsection{Shepp $p$-infinite divisibility for $p$ an even integer}\label{sec:2infdiv}

We have seen that Shepp $p$-infinitely divisible distribution exists for any value of $p \in [1,\infty)$ in Section~\ref{sec:pinfdiv}. When $p$ is an even integer, the subclass of stable distribution $\Sta(\alpha, -1, \lambda, 0)$ parametrized by \eqref{eq:GK} for $\alpha\in (0, 1)$ provides another example. Such a distribution is known to have positive support \cite[Proposition~1.2.11]{Samorodnitsky1994}, and thus there is a well-defined random variable $X_i$ with $1/X_i^p \sim \Sta(\alpha, -1, \lambda, 0)$ for any even integer $p$. The proof is similar to that of Proposition~\ref{prop:exist} and thus omitted.

\begin{proposition}\label{prop:2param}
Let $\alpha_0\in(0, 1)$ and $n \in \mathbb{N}$. Let $X_1, \dots, X_n$ be independent random variables with $1 / X_i^p \sim \Sta(\alpha_0, -1, \lambda_i, 0)$, $X_i>0$, $i=1,\dots,n$, and $T_p\coloneqq S_p(X_1, \dots, X_n)$. Then, for $p \in 2\mathbb{N}$, we have
\[
 \frac{1}{T_p^p} \sim \Sta(\alpha_0, -1, \lambda_1+\dots+\lambda_n, 0).
\]
\end{proposition}
The converse of the above also holds in an appropriate sense:
\begin{corollary}\label{cor:2infdiv}
Let $\alpha\in(0, 1)$ and $p$ be an even integer. If $T_p$ is such that
\[
 \frac{1}{T_p^p} \sim \Sta(\alpha, -1, \lambda, 0), \quad T_p>0,
 \]
then $T_p$ is Shepp $p$-infinitely divisible, i.e., for any $n \in \mathbb{N}$, we have that
 $T_p\deq S_p(X_1, \dots, X_n)$ with independent
 $1 / X_i^p \sim \Sta(\alpha, -1, \lambda/n, 0)$, $X_i>0$, $i=1,\dots,n$.
\end{corollary}

By Proposition~\ref{prop:2param}, we in fact obtain the standard error relation \eqref{eq:S2} for mean zero normal distributions. To show this, we apply the following symmetrization. Denote by $\Ber_s$ the symmetric Bernoulli distribution taking values $\pm 1$ with equal probability.
Let $B_1,\dots, B_n$ and $B$ follow i.i.d.\ $\Ber_s$. One can easily verify that $B \deq B_1 \cdots B_n$,
which enables us to symmetrize $X_i$ by multiplying $B_i$.
Then we get
\[
 S_p(B_1 X_1, \dots, B_n X_n) \deq B \cdot S_p(X_1, \dots, X_n),
\]
leading to a version of Proposition~\ref{prop:2param} for symmetric distributions.

\begin{proposition}\label{prop:2infdivsym}
Let $\alpha_0\in(0, 1)$ and $n \in \mathbb{N}$. Let $1 / X_i^p \sim \Sta(\alpha_0, -1, \lambda_i, 0)$, $X_i>0$, $i =1,\dots,n$, and
 $B_1,\dots,B_n \sim \Ber_s$, with $X_1,\dots,X_n,B_1,\dots,B_n$ independent.
If $T_{B, p} \coloneqq S_p(B_1X_1, \dots, B_nX_n)$, then for $p \in 2\mathbb{N}$,
\[
 \frac{1}{T_{B, p}^p} \sim \Sta(\alpha_0, -1, \lambda_1+\dots+\lambda_n, 0),
 \quad T_{B, p} \deq -T_{B, p}.
\]
\end{proposition}

Now, let $Z_\sigma \sim  \Nor(0, \sigma^2)$.
When $\alpha = 1/2$, $\gamma = 0$, $p = 2$, $\Sta(1/2, -1, \lambda, 0)$ is a L\'evy distribution with scale parameter $\lambda^2$. Hence, $ 1 / Z_{1/\lambda}^2 \sim \Sta(1/2, -1, \lambda, 0)$.
In view of this relation, we have
\[
 S_2(Z_{\sigma_1}, \dots, Z_{\sigma_n}) \sim \Nor(0, 1/\sigma_1^2+ \dots +1/\sigma_n^2),
\]
which is the main result in \cite{Quine1994}. In addition, we can immediately deduce that $Z_{\sigma}$ is Shepp $2$-infinitely divisible since $Z_\sigma \deq S_2(Z_{n\sigma,1}, \dots, Z_{n\sigma,n})$ for any $n \in \mathbb{N}$ and $\sigma>0$. Here $Z_{n\sigma,1}, \dots, Z_{n\sigma,n}$ are $n$ i.i.d.\ copies of $Z_{n\sigma}$.

\subsection{Shepp $\infty$-infinite divisibility}\label{sec:infinfdiv}

We now present a few interesting Shepp $\infty$-infinitely divisible random variables, and show that they are closely related to extreme value distributions.

\begin{definition}[Extreme value distributions]\label{def:extreme}
    We say a random variable follows the \emph{first} or \emph{second} type of extreme value distributions if its distribution function is of the form
    \begin{equation}\label{eq:frechet}
        F_1(x) =
        \begin{cases}
        \exp(-\lambda x^{-r}) & \text{if } x > 0,\\
        0 & \text{if } x \le 0,
        \end{cases}
\quad\text{or}\quad
        F_2(x) =
        \begin{cases}
        1 & \text{if } x> 0,\\
        \exp\bigl(-\lambda(-x)^{r}\bigr) & \text{if } x\le 0,
        \end{cases}
\end{equation}
respectively, for some $r>0$ and $\lambda>0$.    Equivalently, let $U$ and $V$ be random variables with $U^r \sim \Exp(\lambda)$, $U>0$ and $1 / (-V)^r \sim \Exp(\lambda)$, $V<0$. It is easy to verify that $1/U$ and $1/V$ respectively follow the first and second type of extreme value distributions with parameters $(r, \lambda)$. 
\end{definition}

The two types of extreme value distributions are commonly referred to as Fr\'echet distribution \cite{Frechet1927} and reverse Weibull distribution \cite{Weibull1951} in the literature. We favor the $1/U$ and $1/V$ characterizations as $U$ and $V$ are directly relevant to Shepp $\infty$-infinite divisibility, as we will see below.

\begin{proposition}\label{prop:infparam}
Let $n \in \mathbb{N}$, $r>0$, and $\lambda_i>0$, $i = 1, \dots, n$. Let $X_1^r,\dots,X_n^r \sim \Exp(\lambda_i)$, $X_i>0$, be independent. If $T_\infty \coloneqq S_\infty(X_1, \dots, X_n)$, then
\[
 T_\infty^r \sim \Exp(\lambda_1 + \dots+\lambda_n).
\]
\end{proposition}
\begin{proof}
For $x>0$,
\[
 \mathbb{P}\bigl(\max(1/X_1, \dots, 1/X_n)\le x\bigr) = \prod_{i = 1}^n \exp(-\lambda_i x^{-r}) = \exp\biggl(- \sum_{i = 1}^n \lambda_i x^{-r}\biggr).
\]
Hence $\bigl(1 / \max(1/X_1, \dots, 1/X_n)\bigr)^r \sim \Exp(\lambda_1 + \dots+\lambda_n)$.
\end{proof}

As in the $p < \infty$ case, we have a converse in an appropriate sense.
\begin{corollary}
Let $r>0$ and $\lambda>0$. Then $T_\infty$ is such that $T_\infty^r \sim \Exp(\lambda)$ is Shepp $\infty$-infinitely divisible, i.e., for any $n \in \mathbb{N}$, we have $T_\infty \deq S_\infty(X_1, \dots, X_n)$ with independent $X_1^r ,\dots, X_n^r \sim \Exp(\lambda/n)$, $X_i>0$, $i =1,\dots,n$.
\end{corollary}

Likewise, the random variable $V$ in Definition~\ref{def:extreme} also exhibits  Shepp $\infty$-infinite divisibility.

\begin{proposition}
Let $n \in \mathbb{N}$, $r>0$, and $\lambda_i>0$, $i = 1, \dots, n$. Let  $X_1,\dots,X_n$ be independent with $1 / (-X_i)^r \sim \Exp(\lambda_i)$ and $X_i < 0$, $i=1,\dots,n$. If $R_\infty \coloneqq S_\infty(X_1, \dots, X_n)$, then
\[
 1 / (-R_\infty)^r \sim \Exp(\lambda_1 + \dots+\lambda_n).
\]
Furthermore, for $r>0$ and $\lambda>0$, if $1 / (-R_\infty)^r \sim \Exp(\lambda)$, then $R_\infty$ is Shepp $\infty$-infinitely divisible, i.e., for any $n\in\mathbb{N}$, we have  $R_\infty \deq S_\infty(X_1, \dots, X_n)$ with
$1 / (-X_i)^r \oversim{iid} \Exp(\lambda/n)$, $X_i<0$.
\end{proposition}

\section{Characterizing distributions with Shepp $p$-product}\label{sec:char}

The stable and max-stable families discussed in Section \ref{sec:infdiv} motivate further study of the distributional fixed-point equation \eqref{eq:distri} --- for what other  distribution does the property \eqref{eq:distri} uniquely characterize it?
For each value of $p$, we consider i.i.d.\ random variables $X_1$ and $X_2$ drawn from an appropriate Shepp $p$-infinitely divisible distribution. We will show that \eqref{eq:distri} holds for some constant $c>0$, and, more importantly, \eqref{eq:distri} characterizes that distribution uniquely. With one small caveat, ubiquitous distributions such as Cauchy, exponential, Laplace, normal may all be characterized in this manner --- if so desired, one may define these distributions as the unique solutions to \eqref{eq:distri} for various values of $p$. 
The small caveat is that we need to introduce an additional parameter, called the regular variation index.

So here we focus primarily on the converse of Section~\ref{sec:infdiv}: To what extent does $S_p(X_1, X_2)$ characterize the distribution of $X_1$? Let $T_p \coloneqq S_p(X_1, X_2)$ in this section for simplicity. To be precise, we study two specific questions:
\begin{enumerate}[\upshape (i)]
\item\label{it:q1} Given the distribution of $T_p$, can we determine the distribution of $X_i$?
\item\label{it:q2} Can we determine the distribution of $X_i$ by solving distributional equation \eqref{eq:distri}, i.e., $cT_p \deq X_1$?
\end{enumerate}

Question~\eqref{it:q1} is easy to answer. For $p = 2$, $T_2$ has been proved to determine the distribution of $X_1$ up to a reflection about zero \cite[Theorem~2.1]{Hamedani2007}, and the argument therein works for all even integers $p$; in other words, $X_i$ can only be partly determined for even integer $p$. For $p=1$, the distribution of $T_1$ uniquely determines the distribution of $X_1$, provided that $\mathbb{P}(X_1=0)=0$ and the characteristic function of $1/X_1$ is nonvanishing. Indeed, $1/T_1=1/X_1+1/X_2$, so the characteristic function of $1/T_1$ is the square of that of $1/X_1$. Similarly, for $p=\infty$, the distribution of $T_\infty$ uniquely determines that of $X_1$, as can be seen directly from their distribution functions.

Question~\eqref{it:q2} is more difficult to answer and represents our effort in this section. For $p = 2$, \cite[Theorem~2.3]{Galambos2003} proved that under some regularity conditions, \eqref{eq:distri} has a solution if and only if $c = 2$ and $X_i\sim \Nor(0, \sigma^2)$ for some $\sigma>0$. We will extend their study to all even integral values of $p$ as well as $p = 1, \infty$. We will see that the \emph{regular variation index} \cite[Chapter~1]{Seneta1976} below plays an important role.

\begin{definition}[Regular variation]\label{def:regvar}
A function $R(x)$ is said to be \emph{regularly varying} as $x\to 0$ if it is real-valued, positive and measurable on a neighborhood of $0$, and if for all $\lambda>0$,
\[
\lim_{x \to 0} \frac{R(\lambda x)}{R(x)}=\lambda^{r}
\]
for some $r\in \mathbb{R}$, called the index of \emph{regular variation}. If necessary, one may replace the constant $0$ with $0^+$ or $0^-$, i.e., taking the right or left limit for $R(\lambda x)/R(x)$.
\end{definition}

As an overview, we summarize our results in this section in Table~\ref{tab:regvar} for different $p$ and $r$, the regular variation index of the distribution function of $X$ at zero. Here $X$ denotes a random variable with the same distribution as $X_1$, $X_2$ in \eqref{eq:distri}. The only case  known before our work is when $p = 2$ and $r = 1$ \cite{Galambos2003}.

\begin{table}[h!]
\centering
\begin{tabular}{c|c|r@{\;$\sim$\;}l}
$p$                & $r$                            & \multicolumn{2}{c}{characterization} \\\hline
$1$                & $1$                           & $X$           & $\Cau(\gamma, \lambda)$ \\
$1$                & $(0,1)\cup(1, 2)$       & $1 / X$     & $\Sta(r, 0, \lambda, 0)$\\
$1$                & $(2, \infty)$               & $1/X$        & $\Nor(0, \sigma^2)$ \\
$2$                & $1$                           & $X$           & $\Nor(0, \sigma^2)$ \\
$2$                & $(0,1)\cup(1, 2)$       & $1/X^2$   & $\Sta(r/2, -1, \lambda, 0)$, $X\deq -X$ \\
$2k$, $k>1$  & $(0, 2k)$                   & $1/X^{2k}$ & $\Sta(r/2k, -1, \lambda, 0)$, $X\deq -X$ \\
$\infty$          & $1$                           & $X$            & $\Exp(\lambda)$\\
$\infty$          & $(0,1)\cup(1,\infty)$ & $X^r$   & $\Exp(\lambda)$, $X>0$
\end{tabular}
\bigskip
\caption{Summary of characterization results; $r$ denotes regular variation index of $X$ at zero.}
\label{tab:regvar}
\end{table}

This is how one should read the table: Using the first row as illustration, assuming standard regularity conditions, when $r = 1$, a necessary and sufficient condition for \eqref{eq:distri} to have a solution is that $X_1$ and $X_2$ independently follow a Cauchy distribution with location $\gamma\in\mathbb{R}$ and scale $\lambda>0$.

We stress that the condition $r=1$ holds for a wide class of distributions: As long as the density is continuous and nonzero at zero, $r = 1$ holds by L'Hospital's rule. In particular, it holds for normal, Laplace, and uniform distributions.

In the following, we formally state our characterization results in Table~\ref{tab:regvar} for $p = 1$, $\infty$ and even integers. We remind readers that we require $c>0$ in the fixed-point equation \eqref{eq:distri}. The proofs can be found in Appendix~\ref{appx:char}.

For $p = 1$, we focus on the most straightforward scenario which captures Cauchy distribution $\Cau(\gamma, \lambda)$; other cases in Table~\ref{tab:regvar} follow a similar proof.

\begin{theorem}[Characterization of Shepp $1$-product]\label{thm:1char}
Let $X_1, X_2$ be independent random variables with a common distribution function $F$ that satisfies
\begin{enumerate}[\upshape (i)]
\item\label{it:1char1} $F(x) - F(-x)$ is regularly varying as $x\to 0^+$ with index $r = 1$;
\item\label{it:1char2} $\lim_{x\to 0^+}\bigl(F(0) - F(-x)\bigr)/\bigl(F(x) - F(0)\bigr) = 1$.
\end{enumerate}
Then the distributional equation \eqref{eq:distri} with $p = 1$ has a solution if and only if $c = 2$ and $X_1\sim \Cau(\gamma, \lambda)$ for some $\gamma\in\mathbb{R}$ and $\lambda>0$.
\end{theorem}

Here we solve \eqref{eq:distri} for even integer $p$ and general $r$. In the following, $\Sta(\alpha, \beta, \lambda, \gamma)$ denotes the stable distribution parametrized by \eqref{eq:GK}. 
\begin{theorem}[Characterization of Shepp $p$-product for even $p$]\label{thm:2char}
Let $p\in 2\mathbb{N}$ and $X_1, X_2$ be independent random variables with a common distribution function $F$ that satisfies 
\begin{enumerate}[\upshape (i)]
\item\label{it:2char1} $F(x) - F(-x)$ is regularly varying as $x\to 0^+$ with index $r\in (0, p)$;
\item\label{it:2char2} $0<F(0)<1$.
\end{enumerate}
Then the distributional equation \eqref{eq:distri} with $p$ has a solution if and only if $c=2^{1/r}$, $1 / X_1^p\sim \Sta(r/p, -1,\lambda, 0)$, and $X_1\deq -X_1$ for some $\lambda>0$.
\end{theorem}

The characterization \cite[Theorem~2.3]{Galambos2003} is a special case of Theorem~\ref{thm:2char} when $p = 2$, $r = 1$: Observe that if $1 / X_1^2\sim \Sta(1/2, -1,\lambda, 0)$ and $X_1\deq -X_1$, then $X_1\sim \Nor(0, \lambda^{-2})$.

Finally, we solve \eqref{eq:distri} for $p=\infty$ and general $r$.
\begin{theorem}[Characterization of Shepp $\infty$-product]\label{thm:infchar}
Let $X_1, X_2$ be independent random variables with a common distribution function $F$ that satisfies 
\begin{enumerate}[\upshape (i)]
\item\label{it:infchar1} $F(x) - F(0)$ is regularly varying as $x\to0^+$ with index $r\in(0, \infty)$;
\item\label{it:infchar2} $F(x)>F(0)$ for $x>0$.
\end{enumerate}
Then the distributional equation \eqref{eq:distri} with $p = \infty$ has a solution if and only if $c = 2^{1/r}$, $X_1>0$, and $X_1^r \sim \Exp(\lambda)$ for some $\lambda>0$.
\end{theorem}

We state a useful special case as corollary, a consequence of the fact that a continuous density that is nonzero at zero would have a distribution  that satisfies the conditions in Theorem~\ref{thm:infchar}.
\begin{corollary}\label{cor:expchar}
Let $X_1, X_2$ be independent random variables with a density function $f$ that is continuous at zero and $f(0)>0$. Then the distributional equation \eqref{eq:distri} with $p = \infty$ has a solution if and only if $c = 2$ and $X_1\sim \Exp(\lambda)$ for some $\lambda>0$.
\end{corollary}

\section{Central limit theorems for Shepp $p$-product}\label{sec:clt}

We will establish central limit theorems for Shepp $p$-product when $p = 1$, $\infty$, or is an even integer. In other words,  $n^d S_p(X_1,\dots,X_n)$ converges to some limiting distribution regardless of the original distribution of the random variables $X_1,X_2, \dots,$ which are only required to be i.i.d. Here both the power $d$ and the limiting distribution depend only on $p$. In particular, we will establish the following
limits and rates of convergence:
\begin{enumerate}[\upshape (i)]
\item\label{it:1} $p = 1$, $nS_1(X_1, \dots, X_n)\dto\Cau(\gamma, \lambda)$ at rate $O(n^{-1} \log^2 n)$;
\item\label{it:2} $p = 2$, $n S_2(X_1, \dots, X_n) \dto \Nor(0, \sigma^2)$ at rate $O(n^{-\delta})$, $\delta\in (1/2, 1)$;
\item\label{it:infty} $p = \infty$, $n S_\infty(X_1,\dots,X_n)\dto \Exp(\lambda)$ at rate $O(n^{-\delta})$, $\delta\in (0, 1]$;
\end{enumerate}
under mild standard assumptions on regularity and the regular variation index $r$. Notably, $\delta$ only depends on the underlying distribution and the value of $p$. Our actual stated results (see Theorems~\ref{thm:1clt}--\ref{thm:infrate}) will go further than \eqref{it:1}, \eqref{it:2}, \eqref{it:infty}, which represent only selected special cases.

To avoid overly lengthy theorem statements, we first discuss the central limit theorems for various Shepp products, deferring their rates of convergence to the latter half of this section. All proofs may be found in Appendix~\ref{appx:clt}.

The central limit theorem for Shepp $1$-product says that $nS_1(X_1, \dots, X_n)$ converges to a Cauchy law. While we are stating the next result for $r=1$ and the Cauchy distribution, by adjusting the value of $r$, we may prove similar results where the limiting distribution is some other Shepp $1$-infinitely divisible distributions (cf.\ Corollary~\ref{cor:1infdiv}) like Chi-square (cf.\ Corollary~\ref{cor:chisq}).
\begin{theorem}[Central limit theorem for Shepp $1$-product]\label{thm:1clt}
Let $X_1, X_2,\dots$ be independent random variables with a common distribution function $F$ such that $F(x) - F(-x)$ is regularly varying as $x\to 0^+$ with index $r = 1$ and that
\begin{equation}
  \lim_{x\to 0^+}\frac{F(0) - F(-x)}{F(x) - F(0)} = 1.\label{eq:1clt}
\end{equation}
Let
\begin{align}
\tau &\coloneqq
\lim_{n\to\infty}\frac{n}{
    \inf\bigl\{
    x>0:
    \mathbb{P}\bigl(|X_1^{-1}|>x\bigr)
    \le \frac{2}{n\pi}
    \le
    \mathbb{P}\bigl(|X_1^{-1}|\ge x\bigr)
    \bigr\}}\in(0,\infty), \label{eq:1clt3} \\
   \gamma &\coloneqq  \lim_{n \to \infty}\int_{|x|>\tau/n} \frac{1}{x} \, dF(x) <\infty. \label{eq:1clt2}
\end{align}
Then
\[
nS_1(X_1, \dots, X_n)\dto\Cau\Bigl(\frac{\gamma\tau^2}{\gamma^2\tau^2 + 1}, \frac{\tau}{\gamma^2\tau^2 + 1}\Bigr).
\]
\end{theorem}
The condition~\eqref{eq:1clt} regularizes the behavior of $F$ around zero, whereas conditions~\eqref{eq:1clt3} and \eqref{eq:1clt2} ensure that the asymptotic distribution is nondegenerate.

Our next central limit theorem is a generalization of the classical case $p = 2$ \cite{Quine1999} to all even integer values. In such general $p$, the limiting distribution is the $-p$ power of a stable distribution as in Proposition~\ref{prop:2infdivsym}.
\begin{theorem}[Central limit theorem for Shepp $p$-product, $p$ even]\label{thm:2clt}
Let $X_1 , X_2 , \dots$ be independent random variables with a common distribution function $F$ such that $F(0)\in(0, 1)$; $F(x) - F(-x)$ is regularly varying as $x\to 0^+$ with index $r = p\alpha$, $\alpha\in(0, 1)$; and that
\begin{equation}\label{eq:2clt}
 \sigma\coloneqq \lim _{n \to \infty} \frac{n^{1/\alpha}}{\inf \{x: F(x^{-1/p})-F(-x^{-1/p})= 1 / \bigl(n\Gamma(1-\alpha)\cos(\pi\alpha/2)\bigr)\}} \in(0, \infty).
\end{equation}
Then
\[
n^{1/(p\alpha)}S_p(X_1, \dots, X_n)
\dto U
\]
where $1 / U^p \sim \Sta(\alpha, -1, \sigma^{-\alpha}, 0)$ and $U \deq -U$.
\end{theorem}

The result in \cite[Theorem~1]{Quine1999} is a special case of Theorem~\ref{thm:2clt} with $p  =2$ and $\alpha = 1/2$. However, even in this special case, Theorem~\ref{thm:2clt} goes further --- we have replaced the symmetricity assumption in \cite[Theorem~1]{Quine1999} by a much weaker one, namely, $F(0)\in (0, 1)$.

It remains to state the $p = \infty$ case. This does not readily follow from  Theorem~\ref{thm:2clt} via any limiting argument with $p \to \infty$, at least not without additional assumptions that weaken the result.
\begin{theorem}[Central limit theorem for Shepp $\infty$-product]\label{thm:infclt}
Let $X_1, X_2, \dots$ be independent random variables with a common distribution function $F$ such that $F(x)>F(0)$ for all $x>0$; $F(x) - F(0)$ is regularly varying as $x\to0^+$ with index $r>0$; and that
\begin{equation}\label{eq:infclt}
 \lambda \coloneqq  \lim _{n \to \infty} \frac{\inf \{x: F(1 / x)-F(0) \le 1 / n\}}{n^{1/r}} \in(0, \infty).
\end{equation}
Then
\[
n^{1/r} S_\infty(X_1,\dots,X_n)\dto U
\]
where $U>0$ satisfies $U^r\sim \Exp(\lambda^r)$.
\end{theorem}
Note that \(U\) is the random variable in Definition~\ref{def:extreme} with parameters \(r\) and \(\lambda^r\). We state the most important special case below.
\begin{corollary}
Let $X_1, X_2, \dots$ be independent random variables with a common density $f$ continuous at zero and $f(0)>0$.
Then
\[
n f(0) S_\infty(X_1,\dots,X_n)\dto Y
\]
where  $Y \sim \Exp(1)$.
\end{corollary}

We now proceed to provide rates of convergence in Kolmogorov distance for Theorems~\ref{thm:1clt}, \ref{thm:2clt}, and \ref{thm:infclt}. These will depend on a more refined regular variation condition for the distribution function $F$ at zero. It is perhaps worth emphasizing that our results cannot be deduced from Berry--Esseen as is no constant dominance relation between our rates and the classic Berry--Esseen bound $O(1/\sqrt{n})$.

We will need the notion of \emph{domain of normal attraction} from the theory of stable distributions \cite[Section~35]{Gnedenko1954} for finite $p$, but as this has no known analogue for $p = \infty$, we will develop it in our article. It is easiest to illustrate this notion with the Cauchy law: We say the distribution function $G$ belongs to the domain of normal attraction of $\Cau(0,1)$ if for $x>0$,
\[
 G(-x) = \Bigl(\frac{1}{\pi\tau} + a_1(x)\Bigr)
 \frac{1}{x}, \qquad 1 - G(x) = \Bigl(\frac{1}{\pi\tau} + a_2(x)\Bigr)
 \frac{1}{x},
\]
where $\lim_{x\to \infty}a_1(x) = \lim_{x\to \infty}a_2(x) = 0$. Intuitively, the domain of normal attraction restricts the behavior of tails on both sides. An equivalent condition is $\tau S_n/n - a_n\dto \Cau(0, 1)$, with $S_n$ the sum of i.i.d.\ random variables following $G$ \cite[Section~35, Theorem~5]{Gnedenko1954}.

Let $G$ be the distribution function of $1/X_i$ and $F$ the distribution function of $X_i$, and suppose the functions $a_1$ and $a_2$ are continuous on $\mathbb{R}$. Then another equivalent condition for $G$ to belong to the domain of normal attraction of $\Cau(0,1)$ is that for $x>0$,
\begin{equation}\label{eq:1dna}
 F(0) - F(-1/x) = \Bigl(\frac{1}{\pi\tau} + a_1(x)\Bigr)
 \frac{1}{x}, \qquad F(1/x) - F(0) = \Bigl(\frac{1}{\pi\tau} + a_2(x)\Bigr)
 \frac{1}{x}.
\end{equation}
Recall that a function $f:(0,\infty)\to\mathbb{R}$ is
ultimately monotone if there exists $x_*>0$ such that $f$ is
monotone on $[x_*,\infty)$. 
\begin{theorem}[Convergence rate for Shepp $1$-product]\label{thm:1rate}
Let $X_1 , X_2 , \dots$ be independent random variables with a common continuous distribution function $F$ that satisfies \eqref{eq:1dna} for some $\tau\in (0, \infty)$ and all $x>0$.  Let the tail sum and tail difference, i.e.,
\[
S(x) = \frac{a_1(x)+a_2(x)}{x}
\quad\text{and}\quad
D(x) = \frac{a_2(x)-a_1(x)}{x},
\]
be ultimately monotone, with
\[
\int_1 ^{x} y|S(y)|\, d y=O(\log^2x), \qquad|D(x)|=O(x^{-2}\log^2 x)\quad \text{as}\quad x \to \infty.
\]
Then $F_n$, the distribution function of $nS_1(X_1,\dots,X_n)$, satisfies
\[
 \sup_{x\in\mathbb{R}}|F_{n}(x)-\mathbb{P}( C\le x)|
 = O(n^{-1} \log^2 n) \quad \text{ as }\quad n\to \infty,
\]
where $C\sim \Cau\bigl(\gamma\tau^2 / (\gamma^2\tau^2 + 1), \tau / (\gamma^2\tau^2 + 1)\bigr)$ with $\gamma \coloneqq\int_0^\infty D(x)\,dx$.
\end{theorem}

Again, we consider the domain of normal attraction and apply \cite[Section~35, Theorem~5]{Gnedenko1954}. Specifically, we consider the random variables $X_i$ whose common distribution function satisfies
\begin{equation}\label{eq:2dna}
 F(x^{-1/p}) - F(-x^{-1/p}) = \bigl(c_2 \sigma^{-\alpha} + a_1(x)\bigr)x^{-\alpha},
 \qquad c_2 = \frac{1}{\Gamma(1-\alpha)\cos(\pi \alpha/2)}.
\end{equation}
As shown in Appendix~\ref{appx:clt}, Equation \eqref{eq:2dna} is the necessary and sufficient conditions to ensure that $n^{1/p\alpha} S_p(X_1, \dots, X_n)\dto U$, where $U$ is a random variable symmetric around zero and 
$1/U^p \sim \Sta(\alpha, -1, \sigma^{-\alpha}, 0)$. Recall that $\Sta(\alpha, \beta, \lambda, \gamma)$ is defined by \eqref{eq:GK}. Then we give the rate of convergence for distributions in the domain of normal attraction.

\begin{theorem}[Convergence rate for Shepp $p$-product, $p$ even]\label{thm:2rate}
Let $X_1 , X_2 , \dots$ be independent random variables whose common distribution function $F$ is symmetric and continuous, and \eqref{eq:2dna} hold for some $\sigma\in (0, \infty)$ and $\alpha\in(0, 1)$. Let the tail sum, i.e., 
\[
S(x) = a_1(x)x^{-\alpha}
\]
be continuous, with 
\begin{align*}
&x^{\varepsilon}|S(x)| \text{ ultimately nonincreasing for some }\varepsilon > 0, \\
&|S(x)| = O(x^{-s}) \text{ as } x\to\infty \text{ for some } \alpha < s < \min(1, 2\alpha).
\end{align*}
Then $F_n$, the distribution function of $n^{1/p\alpha} S_p(X_1, \dots, X_n)$, satisfies 
\begin{equation}\label{eq:2rate2}
\sup_{x\in\mathbb{R}}|F_{n}(x)-F_0(x)| = O(n^{1 - s/\alpha}) \quad \text{ as } \quad n\to \infty,
\end{equation}
where $F_0$ is the distribution function of a random variable $U$ symmetric around zero and satisfying $1/U^p \sim \Sta(\alpha, -1, \sigma^{-\alpha}, 0)$.
\end{theorem}
Consider a special case $p = 2$, $\alpha = 1/2$, we have the next corollary, which gives the rate of convergence for Shepp $2$-product to a normal distribution.

\begin{corollary}\label{cor:2rate}
Let $X_1 , X_2 , \dots$ be independent random variables with a common distribution function $F$ that is symmetric and continuous, and for $\sigma\in (0, \infty)$,
\[
 F(1/\sqrt{x}) - F(-1/\sqrt{x}) = \bigl(\sqrt{2/\sigma^2\pi} + a_1(x)\bigr)x^{-1 / 2}.
\]
Let $S(x) = a_1(x)x^{-1 / 2}$ be continuous, with 
\begin{align*}
&x^{\varepsilon}|S(x)| \text{ ultimately nonincreasing for some }\varepsilon > 0, \\
&|S(x)| = O(x^{-s}) \text{ as } x\to\infty \text{ for some } 1/2 < s < 1.
\end{align*}
Then $F_n$, the distribution function of $n S_2(X_1, \dots, X_n)$, satisfies 
\[
\sup_{x\in\mathbb{R}}|F_{n}(x)-\mathbb{P}(Z_\sigma\le x)| = O(n^{1 - 2s})\quad \text{ as } \quad n\to\infty,
\]
where $Z_\sigma \sim \Nor(0, \sigma^2)$.
\end{corollary}
In the following, we establish the $O(1 / \sqrt{n})$ rate of convergence for Shepp $2$-product of non-identical independent random variables. Our bound comes from a finite-sample bound for stable distributions in Wasserstein distance \cite{Johnson2005}, which is attached to Appendix~\ref{appx:clt}.
\begin{theorem}\label{thm:2rate2}
Let $X_1, X_2, \dots$ be a sequence of independent and symmetric random variables with the distribution function $F_i$ such that for some common $\sigma>0$
\[
    F_i(1/\sqrt{x}) - F_i(-1/\sqrt{x}) = \bigl(\sqrt{2/\sigma^2\pi} + a_i(x)
    \bigr)x^{-\frac{1}{2}}
\]
where $|a_i(x)|\le C_i/x^{\gamma_i}$ holds for any $x\in(0, \infty)$, for some constant $C_i > 0$, $\gamma_i>0$.
If $\sup\{C_i \colon i\in\mathbb{N}_+\} < \infty$ and $\inf\{\gamma_i \colon i\in\mathbb{N}_+\} > \frac{1}{2}$,
we have
\[
n S_2(X_1, \dots, X_n)\dto Z_\sigma, \quad Z_\sigma \sim \Nor(0, \sigma^2)
\]
and $F_n$, the distribution function of $n S_2(X_1, \dots, X_n)$, satisfies 
\[
    \sup_{x\in\mathbb{R}}|F_{n}(x)-\mathbb{P}(Z_\sigma \le x)|
    = O\bigl(\frac{1}{\sqrt{n}}\bigr)\quad \text{ as } \quad n\to\infty.
\]
\end{theorem}

The domain of normal attraction has so far played an important role in determining the rate of convergence to a limiting stable law. 
For extreme-value distributions, the standard domain-of-attraction framework allows general normalizing sequences. Here we specialize this framework to distributions admitting a prescribed normalization of order $n^{1/r}$.
In this section, first focusing on the convergence of sample maximum, we propose the domain of normal attraction of an extreme value distribution and derive the convergence rates for distributions belonging to such a domain of normal attraction. Secondly, we apply our results to give the rate of convergence for Shepp $\infty$-product.

Our next definition is a maximum analogue of \cite[Section~35, Equation~21]{Gnedenko1954}. For simplicity of the presentation, we focus on the first type of extreme distribution.
\begin{definition}[Domain of normal attraction]\label{def:dna}
Let $X_1, X_2, \dots$ be independent random variables with a common distribution function $G$. We say $G(x)$ belongs to the \emph{domain of normal attraction of an extreme value distribution}, if for some $\lambda, r>0$ and a sequence $a_n$, $n\ge 1$, the following relation holds:
\[
\frac{\max(X_1, \dots, X_n)}{n^{1/r}} - a_n \dto \frac{1}{U},
\]
where $U > 0$ and $U^r \sim \Exp(\lambda^r)$, following the first type of extreme value distribution.
\end{definition}

The next proposition gives a sufficient condition for $G$ to belong to the domain of normal attraction defined above.
\begin{proposition}\label{prop:infdna}
Let $X_1, X_2, \dots$ be independent random variables with a common distribution function $G$, and assume that $G$ satisfies
\begin{equation}\label{eq:infdna}
 1 - G(x) = \bigl(\lambda^r + a_1(x)\bigr)x^{-r},\quad x>0
\end{equation}
where $\lambda>0$ is a constant and $a_1(x)$ is a function such that $\lim_{x \to \infty}a_1(x) = 0$. Then $G$ is in the domain of normal attraction of an extreme value distribution with $a_n = 0$ and $U>0$, $U^r \sim \Exp(\lambda^r)$.
\end{proposition}
In fact, there have been many works on the rate of convergence to $1/U$, the first type of extreme value distribution; we refer readers to \cite{Reiss1989, Falk1986, Smith1982, Galambos1978} for previous results. Our next theorem gives a new rate of convergence for distributions belonging to the domain of normal attraction.
\begin{theorem}\label{thm:infrate}
Let $X_1, X_2, \dots$ be independent random variables with a common distribution function $G$ satisfying \eqref{eq:infdna} with $r>0$. Let $a_1(x)$ in \eqref{eq:infdna} be ultimately nonnegative and satisfy $|a_1(x)|\le C_0/|x|^\eta$ for any $x>0$ and some positive constants $C_0$, $\eta$. Then, we have
\[
 \sup_{x\in\mathbb{R}}\Bigl|\mathbb{P}\Bigl(\frac{\max_i(X_1, \dots, X_n)}{ n^{1/r}}\le x\Bigr) - F_{r}(x)\Bigr| \le \frac{C}{n^s}\quad \text{ as } \quad n\to \infty,
\]
where $F_{r}(x)$ is the distribution function of $1/U$ defined in Proposition~\ref{prop:infdna}, $C>0$ is a constant independent of $n$, and $s = \min(\eta/r, 1)$.
\end{theorem}

As a direct result of Proposition~\ref{prop:infdna} and Theorem~\ref{thm:infrate}, we establish the rate of convergence for Shepp $\infty$-product. Its proof is omitted for simplicity.
\begin{theorem}[Convergence rate for Shepp $\infty$-product]
Let $X_1, X_2, \dots$ be independent random variables with a common continuous distribution function $F$, $\lambda$,$r$ be positive constants, and $U$ be a positive random variable such that $U^r \sim \Exp(\lambda^r)$. Assume that
\[
F(1/x) - F(0) = \bigl(\lambda^r + a_1(x)\bigr)x^{-r}, \quad x>0
\]
where $a_1(x)$ is ultimately nonnegative and satisfies $|a_1(x)|\le C_0/|x|^\eta$ for positive constants $C_0, \eta$. We have
\[
\sup_{x\in\mathbb{R}}\bigl|\mathbb{P}\bigl(n^{1/r} S_\infty(X_1, \dots, X_n)\le x\bigr) - \mathbb{P}(U \le x)\bigr| \le \frac{C}{n^s}
\]
with some $C>0$ and $s = \min(\eta/r, 1)$.
\end{theorem}

\section{Tests based on Shepp $p$-product characterizations}\label{sec:test}

The characterization results in Section~\ref{sec:char} provide a natural basis for testing certain distributions from i.i.d.\ observations: if the underlying distribution is a solution to \eqref{eq:distri}, we expect the empirical distributions of both sides of \eqref{eq:distri} are similar; otherwise they should be distinguishable from each other. In this section, we construct the integral-type and the Kolmogorov-type statistics to measure the discrepancy between two empirical distribution functions. While this procedure is valid for any solution distribution to \eqref{eq:distri}, we mainly elaborate the testing procedure for the family of exponential distributions and figure out their asymptotic distributions under the null. It is worth noting that the proposed test statistics are scale-invariant, which releases us from a first-stage scale parameter estimation (i.e., data standardization) in the classic Kolmogorov--Smirnov test. We will also propose Cauchy tests similarly at the end of this section.

To analyze the power under alternative distributions, we adopt the notion of local Bahadur efficiency \cite{Nikitin1995} as our criterion, which, intuitively speaking, measures the efficiency to reject the null under an alternative hypothesis. We compute local Bahadur efficiency for our exponentiality tests. In addition, an extensive simulation is carried out against various alternatives. Both theory and simulation confirm the superiority of the proposed $\mathit{IU}_n^E$ statistic over Kolmogorov--Smirnov test.
\subsection{Exponentiality tests and their asymptotic distributions}\label{sec:exp1}
Define $\mathscr{F}_{\mathscr{E}}$ to be the family of distribution functions such that $F(x)>F(0)$ for $x>0$, and that $F(x) - F(0)$ is regularly varying as $x\to0^+$ with index $1$. For elements in $\mathscr{F}_{\mathscr{E}}$,
\[
\frac{2}{\max(1/X_1, 1/X_2) } \deq X_1
\]
holds only for the distribution functions of exponential distributions $\Exp(\lambda)$ with $\lambda>0$ due to Theorem~\ref{thm:infchar} and Corollary~\ref{cor:expchar}.

Let $X_i$ be i.i.d.\ samples from a distribution function $F$, and $F_n$ be the respective empirical distribution function. Define the composite hypothesis $H_0: F \sim \Exp(\lambda)$ against the alternative hypothesis $H_1: F \in \mathscr{F}_{\mathscr{E}}$ but $H_0$ is false. For brevity, let $k_E(x, y) = 2/\max(1/x,1/y)$. For any $t \in \mathbb{R}$, we define
\[
V_{n}^E(t)=n^{-2} \sum_{i, j=1}^{n} \mathbb{I}_{\{k_E (X_{i}, X_j)<t\}}, \quad
U_{n}^E(t)= \frac{2}{n(n-1)} \sum_{1 \le i<j \le n} \mathbb{I}_{\{k_E(X_{i}, X_j)<t\}},
\]
and construct two statistics for testing the exponential distributions, i.e.,
\[
\mathit{IU}_n^E = \int_{\mathbb{R}} \bigl(V_n^E(t) - F_n(t)\bigr) \, dF_n(t), \quad
\mathit{DU}_{n}^E=\sup _{t \in \mathbb{R}}|U_{n}^E(t)-F_{n}(t)|.
\]
Then we establish the following result, whose proof can be found in Appendix~\ref{appx:test}. 

\begin{theorem}\label{thm:null}
Under $H_0$, as $n\to\infty$,
\[
    \sqrt{\frac{420n}{11}}\mathit{IU}_n^E \dto \Nor(0, 1), \quad
    \sqrt{n}\mathit{DU}_n^E \dto \sup_t |\eta(t)|.
\]
$\eta(t)$ is a centered Gaussian process whose covariance is given in our proof in Appendix~\ref{appx:test}. 
\end{theorem}
The corresponding tests can then be constructed using rejection rules derived from these asymptotic distributions.
Although the asymptotic distribution of $\sqrt{n}\mathit{DU}_n^E$ is intractable, we could obtain critical values for the test following the procedure in \cite{Volkova2009}.

\subsection{Local Bahadur efficiency and empirical performance of exponentiality tests}
Local Bahadur efficiency is widely-used to evaluate the type-II error: Given a specific alternative, the larger it is, the smaller type-II error we obtain. Its formal definition can be found in \cite{Nikitin1995}.

Following the same procedure in \cite[Section~4]{Nikitin2017}, we derive the local Bahadur efficiency in Table~\ref{tab:eff} for the proposed exponentiality tests, considering Weibull and gamma distribution as alternatives. From Table~\ref{tab:eff}, $\mathit{IU}_n^E$ achieves a relatively high efficiency.
\begin{table}[ht]
\centering
\begin{tabular}{c|cc|cc}
\hline
Test & \multicolumn{2}{c|}{$\mathit{IU_n^E}$} & \multicolumn{2}{c}{$\mathit{DU_n^E}$} \\
\hline
Alternative & Weibull & Gamma & Weibull& Gamma \\
\hline
Efficiency & 0.6971 & 0.7901 & 0.1581 & 0.1737 \\
\hline
\end{tabular}
\bigskip
\caption{Local Bahadur efficiency of exponentiality tests.}
\label{tab:eff}
\end{table}

To support our theoretical finding, we carry out an empirical power study and include Kolmogorov--Smirnov test as a baseline. Since $\lambda$ is assumed unknown in the data, we use the MLE estimator $n/\sum_i X_i$ Kolmogorov--Smirnov test in place of $\lambda$. We experiment $\Exp(1)$, $\Exp(5)$ as null hypotheses. As for the alternative, we consider a broad class of distributions adapted from \cite{Ascher1990}, including gamma distribution with shape parameter $a$ and scale parameter 1 (denoted by $\Gam(a,1)$), Weibull distribution ($\Wb(a)$) with density
\[
f(x;a) = a x^{a-1} \exp(-x^a),
\]
beta distribution with shape parameters $a,b$ ($\Bet(a, b)$), and lognormal distribution with scale parameter $a$ and location parameter zero ($\LN(a)$). To be specific, we experiment three distributions with monotonically decreasing hazard rates($\Gam(0.5,1)$, $\Gam(0.7,1)$, $\Wb(0.8)$), eight with monotonically increasing hazard rates ($\Uni(0,1)$, $\Gam(1.5,1)$, $\Gam(2,1)$, $\Gam(4,1)$, $\Wb(1.2)$, $\Wb(1.5)$, $\Bet(1,2)$, $\Bet(2,1)$), and three whose hazard rates are non-monotonic ($\LN(0.6)$, $\LN(1)$, $\Bet(0.5, 1)$). In addition, we test contaminated alternatives, namely a.) 5 samples from $\Exp(1/5)$ and 45 samples from $\Exp(1)$ and b.) 5 samples from $\Exp(1/3)$ and 45 samples from $\Exp(1)$.

The sample size is fixed at 50 and 1000 values of each test statistic are simulated for each specified distribution. A type-I error of $0.05$ is utilized, and rejection rates in both null and alternative hypotheses are presented in Table~\ref{tab:exp}.

Simulation results support the efficiency given in Table~\ref{tab:eff}. In addition, one could see from Table~\ref{tab:exp} that the proposed integral-type statistic $\mathit{IU}_n^E$ has good control of type-I error and remarkable power under various alternatives. 
By contrast, Kolmogorov-type statistics can have relatively low efficiency \cite{Nikitin1995, Nikitin2010_largedeviation}, and the corresponding tests in Table~\ref{tab:exp} exhibit comparatively modest empirical power for several alternatives.

\begin{table}[h]
\centering
\begin{tabular}{lllll}
Types & Distributions & $\mathit{IU}^E_n$ & $\mathit{DU}^E_n$ & KS \\
\hline
\hline
\multirow{2}{10em}{Null hypotheses} & $\Exp(1)$ & 0.063 & 0.014 & 0.009 \\
&$\Exp(5)$ & 0.063 & 0.012 & 0.007 \\
\hline
\multirow{3}{15em}{Monotonic decreasing hazard} & $\Gam(0.5,1)$ & 0.917 & 0.015 & 0.629 \\
&$\Gam(0.7,1)$ & 0.377 &  0.005 & 0.117 \\
& $\Wb(0.8)$ & 0.332& 0.003& 0.142 \\
\hline
\multirow{8}{15em}{Monotonic increasing hazard} & $\Uni(0, 1)$ & 0.837 & 0.755 & 0.665 \\
&$\Gam(1.5,1)$ & 0.547 & 0.150 & 0.110 \\
&$\Gam(2,1)$ & 0.932 & 0.448 & 0.472 \\
&$\Gam(4,1)$ & 1.000 &  0.996 & 1.000 \\
&$\Wb(1.2)$ &0.326& 0.084& 0.049 \\
&$\Wb(1.5)$ &0.904& 0.379& 0.442 \\
&$\Bet(1,2)$ & 0.349 & 0.108 & 0.152 \\
&$\Bet(2,1)$ & 1.000 & 1.000 & 1.000 \\
\hline
\multirow{3}{15em}{Non-monotonic hazard} &
$\Bet(0.5,1)$ & 0.466 & 0.020 & 0.094 \\
&$\LN(0.6)$ & 1.000 & 0.962& 0.985 \\
&$\LN(1)$ & 0.394 & 0.120 & 0.058 \\
\hline
\multirow{2}{15em}{Contaminated exponential} & 5 $\Exp(1/5)$ and 45 $\Exp(1)$ & 0.066& 0.012& 0.136 \\
& 5 $\Exp(1/3)$ and 45 $\Exp(1)$ & 0.064& 0.011& 0.019 \\
\hline
\hline
\end{tabular}
\bigskip
\caption{Rejection rates of exponentiality tests using $\mathit{IU}^E_n$, $\mathit{DU}^E_n$, and KS. KS denotes Kolmogorov--Smirnov test.}
\label{tab:exp}
\end{table}

\subsection{Tests for the Cauchy distribution}

Define $\mathscr{F}_{\Cau}$ to be the family of distribution functions such that conditions \eqref{it:1char1}, \eqref{it:1char2} in Theorem~\ref{thm:1char} are satisfied. For elements in $\mathscr{F}_{\Cau}$,
\[
\frac{2X_1X_2}{X_1 + X_2} \deq X_1
\]
holds only for the distribution functions of Cauchy distributions
$\Cau(\gamma, \lambda)$ with $\gamma\in\mathbb{R}$, $\lambda>0$ due to Theorem~\ref{thm:1char}.

Similar to the exponential case, we work on the composite null hypothesis $H_0: F \sim \Cau(0, \lambda)$
against the alternative hypothesis $H_1: F \in \mathscr{F}_{\Cau}$ but $H_0$ is false.
Let $k_C(x, y) = 2xy/(x+y)$ and define for any $t \in \mathbb{R}$,
\[
V_{n}^C(t) \coloneqq n^{-2} \sum_{i, j=1}^{n} \mathbb{I}_{\{k_C (X_{i}, X_j)<t\}}, \quad 
U_{n}^C(t) \coloneqq \frac{2}{n(n-1)} \sum_{1 \le i<j \le n} \mathbb{I}_{\{k_C(X_{i}, X_j)<t\}}.
\]
Next, we construct two statistics for testing the Cauchy distribution
\[
    \mathit{IU}_n^C \coloneqq \int_{\mathbb{R}} \bigl(V_n^C(t) - F_n(t)\bigr) \, dF_n(t), \quad
    \mathit{DU}_{n}^C \coloneqq \sup _{t \in \mathbb{R}}|U_{n}^C(t)-F_{n}(t)|.
\]
Assume $X_i \oversim{iid} \Cau(0, 1)$ without loss of generality. Using the identical argument as in the proof of Theorem~\ref{thm:null}, we get
\[
    \sqrt{n}\mathit{IU}_n^C \dto \Nor(0, 9\sigma^2), \quad
    \sqrt{n}\mathit{DU}_n^C \dto \sup |\eta(t)|
\]
where $\sigma^2 \approx 0.0052$ and $\eta(t)$ is a centered Gaussian process as in \cite{Silverman1983}. Therefore, the asymptotic distribution for $\mathit{DU}_n^C$ is intractable, but we could estimate its critical values by simulating $\eta(t)$.

\section{Applications of Shepp $p$-product}\label{sec:app}
In this section, we discuss some applications of the Shepp $p$-product in biology and physics. We briefly propose new methods in two specific setups, but consider a more concrete study as future work.
\subsection{Statistics in alignment-free sequence analysis}\label{sec:align}
The conceptual goal in alignment-free sequence analysis is to test whether two sequences are generated independently or not, without relying on the biological knowledge of any specific patterns in the sequences. The new tests we are going to propose, belong to the category of $k$-tuples (also called $k$-grams or $k$-words) approach. This is a general approach to testing the independence by counting the occurrence of every length-$k$ pattern in a sequence. See \cite{Reinert2009, Song2013, Ren2013} for previous studies. In the following, we utilize the properties of Shepp $p$-product studied so far to construct new test statistics.

First recall the settings described in \cite{Reinert2009}: Suppose that two sequences, $A=A_1A_2\cdots A_{n}$ and $B=B_1B_2\cdots B_{m}$ are composed of letters that are drawn from a finite alphabet $\mathcal{A}$ of size $d$. For $a \in \mathcal{A}$, let $p_{a}$ be the probability of drawing a letter $a$ from $\mathcal{A}$.

The null hypothesis in this problem is that the letters are independent and identically distributed, and that two sequences are independent. If the null hypothesis is true, we expect every word of length $k$, i.e., $w=(w_1 , \dots, w_k)$, has the same frequency in sequences $A$ and $B$. It inspires us to construct test statistics from word-count variables. Define $\bar{n} = n - k + 1$, $\bar{m} = m - k + 1$, and $p_{w}=\prod_{i=1}^k p_{w_i}$ the probability of the occurrence of $w$. The word-count variables (after normalization) are
\[
     X_{w}=\frac{1}{\sqrt{n}}\sum_{i=1}^{\bar{n}} \mathbb{I}_{\{A_i =w_1 , \dots, A_{i+k-1}=w_k \}}-\frac{\bar{n}}{\sqrt{n}} p_{w}, \quad
     Y_{w}=\frac{1}{\sqrt{m}}\sum_{i=1}^{\bar{m}} \mathbb{I}_{\{B_i =w_1 , \dots, B_{i+k-1}=w_k \}}-\frac{\bar{m}}{\sqrt{m}} p_{w}.
\]
Lastly, we construct a set of test statistics,
\[
 D_p^S \coloneqq \sum_{w \in \mathcal{A}^k } S_p(X_{w}, Y_{w})
\]
where the superscript $S$ stands for Shepp. The index $p$ can be $1$, $\infty$, or even numbers. Their asymptotic distributions, formally proved in Appendix~\ref{appx:app}, are given below. Throughout this section, we consider fixed $d, k$ and let $m, n$ go to infinity.
\begin{theorem}\label{thm:align1}
Let $Z_1 =(Z_{1,1}, \dots, Z_{1, d^k -1}) \sim \Nor(0, \Sigma_0)$ and $Z_2 =(Z_{2,1}, \dots, Z_{2, d^k -1}) \sim \Nor(0, \Sigma_0)$ be two independent $(d^k -1)$-dimensional normal vectors with $\Sigma_0$ defined in Appendix~\ref{appx:app}, and for $i = 1,2$,
\[
 Z_{i, d^k}=-\sum_{j = 1}^{d^k - 1} Z_{i, j}. 
\]
Let
\[
 D_p^{\lim} = \sum_{j = 1}^{d^k} S_p(Z_{1, j}, Z_{2, j}).
\]
Under the null hypothesis and mild conditions in the proof, we have $D_p^S\dto D_p^{\lim}$ for $p = 1, \infty$ and even numbers. 
\end{theorem}
While $D_2^{\lim}$ follows a mean zero normal distribution \cite{Reinert2009}, for $p\ne 2$, the distribution of $D_p^{\lim}$ is in general intractable. As an alternative, we suggest computing the critical values of $D_p^{\lim}$ by sampling.

For $p = 1$, we can construct another statistic in the form of Shepp $1$-product, different from $D^S_1$, 
\[
C_1 ^S = \sum_{w, w^\prime \in \mathcal{A}^k} \frac{X_{w}X_{w^\prime}}{X_{w}Y_{w^\prime} + X_{w^\prime}Y_{w}}.
\]
Following the assumptions and the proof of Theorem~\ref{thm:align1}, it is straightforward to conclude that $C_1^S\dto C_1^{\lim}$ as $m,n$ go to infinity, where
\[
C_1^{\lim} = \sum_{i,j = 1}^{d^k} S_1\biggl(\frac{Z_{1,i}}{Z_{2,i}}, \frac{Z_{1,j}}{Z_{2,j}}\biggr).
\]
Again, we suggest computing its critical values by sampling.

To summarize, we propose two kinds of statistics $D_p^S$ and $C_1^S$ for the alignment-free sequence test. Our test statistics are applicable in general.  However, it is unclear about its power in either theory or simulation study, which we consider as future work.

\subsection{Asymptotic conductivity of random resistor network}\label{sec:media}
The random resistor network is a hierarchical model that represents the natural structure underlying many materials in physics. We are interested in the conductance of the network, which can be expressed as a weighted sum of $n$ i.i.d.\ random variables. As a result, central limit theorems for a certain summation of i.i.d.\ random variables have been considered extensively in modeling conductance in random media. Some concrete examples bridging the gap between the convergence theorems and applications in random media are described in \cite[Section~4.2]{Chen2011} and \cite{Li1999}.

Let's consider a simple case of random resistors in series, i.e., random resistors with i.i.d.\ conductance $X_i$ in series. Denote by $L(X_1, \dots, X_n)$ the total conductance for $n$ random resistors. Since resistances add for components placed in series, and resistance and conductance are inverses \cite[Section~4.2]{Chen2011}, we have for a pair of components that $L(X_1, X_2) = (X_1^{-1} + X_2^{-1})^{-1} = S_1(X_1, X_2)$. For $n$ resistors, $L(X_1, \dots, X_n) = S_1(X_1, \dots, X_n)$ holds analogously. For $X_i$ which follows some positive distribution, observe that 
\[
S_1(X_1, \dots, X_n) = \bigl(S_2(\sqrt{X_1}, \dots, \sqrt{X_n})\bigr)^2.
\]
Then, the proof of Theorem \ref{thm:2clt} (given $p = 2$) yields
\[
 n^{1/\alpha}S_1(X_1, \dots, X_n)\dto U, \quad 1/U \sim  \Sta(\alpha, -1, \sigma^{-\alpha}, 0)
\]
under the corresponding regular variation and normalization assumptions, with $\alpha$ and $\sigma$ defined accordingly. That is, the total conductance of a random resistor system converges to the inverse of some positive stable distributions. Specific asymptotic distribution for $L(X_1, \dots, X_n)$ is determined by the regular variation index $\alpha$ around zero.

\section{Discussion}
In this paper, we extensively study the statistical properties of a new type of product, called Shepp $p$-product. We develop a new notion of infinite divisibility with respect to Shepp $p$-product, tackle the characterization problem \eqref{eq:distri}, and derive the asymptotic distributions along with the respective rates of convergence for different Shepp $p$-products of i.i.d.\ random variables.

We propose assorted applications of Shepp $p$-product, including hypothesis testing for distributions, alignment-free sequence analysis, and random resistor network. For alignment-free sequence analysis, a concrete application is to apply our tests to Forster's Covid data set for phylogenetic network analysis \cite{Forster2020}, which is considered as future work.

It is also tempting to generalize the results in Sections \ref{sec:char} and \ref{sec:clt} to dependent random variables. For instance, a much weaker condition for
\[
    c S_1(X_1, X_2) \deq X_1
\]
to have a solution is that $(X_1, X_2)$ follows a joint Cauchy distribution in the sense of \cite{Pillai2016}; it is a direct application of the main theorem in \cite{Pillai2016}. That said, the sufficient part of Theorem~\ref{thm:1char} can be considerably relaxed. Admittedly, it is unclear to what extent those results can be generalized, and we consider it a potential research direction.

\section*{Acknowledgement}
The authors thank Mathias Drton and the other participants in the Algebraic Statistics Online Seminar for their constructive comments on an earlier version of this work.

\appendix

\section{Preliminary Results}
In this section, we provide preliminary results used in our proofs.

Many of our proofs rely on a more involved version of Theorem~\ref{thm:GK}, which we state below for completeness.
\begin{theorem}\label{thm:GKcomplete}
In order that a distribution is stable, it is necessary and sufficient that its characteristic function $\varphi$ be represented by
\[
\log \varphi(t)=i \gamma t+\int_{u \ne 0}\Bigl(e^{i t u}-1-\frac{i t u}{1+u^2 }\Bigr) \, d \Lambda(u)
\]
with $\alpha\in(0, 2]$, $\gamma\in\mathbb{R}$, and
\[
\Lambda(u)=
\begin{cases} c_1/|u|^{\alpha} & \text{if } u<0,\\
 -c_2/ u^{\alpha} & \text{if }  u>0,
\end{cases}
\]
where $c_1 , c_2\ge 0$ and $c_1 +c_2 >0$. Furthermore, the above representation can be simplified to \eqref{eq:GK}: For $\alpha \ne 1$, the parameters $\lambda$ and $\beta$ in \eqref{eq:GK} are deduced by
\[
\lambda=(c_1 +c_2 ) \cos \Bigl(\frac{\pi \alpha}{2}\Bigr)\Gamma(1-\alpha), \quad
\beta=\frac{c_1 -c_2 }{c_1 +c_2 };
\]
for $\alpha = 1$, similarly, $\lambda$ and $\beta$ are deduced by
\[
\lambda=(c_1 +c_2 ) \frac{\pi}{2}, \quad \beta=\frac{c_2 -c_1}{c_1 +c_2}.
\]
\end{theorem}
The quantities $c_1, c_2$ play a crucial role in the central limit theorem below, which is adapted from \cite[Section~35]{Gnedenko1954}.
\begin{theorem}\label{thm:GKclt}
Let $X_1 , X_2 , \dots$ be independent random variables with a common distribution function $F$. In order that 
\[
\frac{1}{b_n}\sum_{i = 1}^n X_i - a_n
\]
converges weakly to a stable distribution with some suitable $a_n$, $b_n$, it is necessary and sufficient that 
\begin{gather}
    \frac{F(-x)}{1-F(x)} \to \frac{c_{1}}{c_{2}} \quad\text{ as } x \to \infty, \label{eq:clt1}\\
    \frac{1-F(x)+F(-x)}{1-F(k x)+F(-k x)} \to k^{\alpha} \quad \text { as } x \to \infty\label{eq:clt2}.
\end{gather}
By choosing
\begin{align*}
    a_n &= 
    \begin{cases}
    0 & 0<\alpha<1, \\
    \frac{n}{b_n}\int_{|x|<b_n} x \, d F(x) & \alpha = 1, \\
    \frac{n}{b_n}\int_{-\infty}^{\infty}x \, d F(x) & 1<\alpha\le 2,
    \end{cases} \\
    b_n &= \inf\Bigl\{x\colon \mathbb{P}(|X_1|>x)\le \frac{c_1+c_2}{n}\le \mathbb{P}(|X_1|\ge x)\Bigr\},
\end{align*}
the limiting stable distribution is uniquely determined by $c_1$, $c_2$, $\alpha$, and $\gamma = 0$.
\end{theorem}

The following lemma gives the uniform convergence of distribution functions, which will be used constantly in the main proofs.
\begin{lemma}\label{lemma:cdf}
Let a sequence of distribution functions $F_n(x)$ converge weakly to a continuous distribution function $F(x)$ and a sequence of real numbers $a_n$ converge to $a$. Then we have
\[
 \lim_{n\to\infty}F_n(a_n) = F(a).
\]
\end{lemma}
\begin{proof}
By triangle inequality, we have
\begin{align*}
 \lim_{n\to\infty}|F_n(a_n) - F(a)|&\le \lim_{n\to\infty}|F_n(a_n) - F(a_n)| + \lim_{n\to\infty}|F(a_n) - F(a)|.
\end{align*}
Since $F(x)$ is continuous, $F_n(x)$ converges to $F(x)$ uniformly (see, for example, \cite[Lemma~2.10.1]{Galambos1978}), i.e., $\lim_{n\to\infty}\sup_{x\in\mathbb{R}}|F_n(x) - F(x)| = 0$. This completes the proof. 
\end{proof}

In some of our proofs, we utilize a new notion of distance, called Wasserstein $\nu$-distance, in order to apply certain results from previous works. Specifically, the Wasserstein $\nu$-distance is defined by
\[
    d_{\nu}(X_0, Y_0) \coloneqq \Bigl(\inf _{(X, Y)} \mathbb{E}|X-Y|^{\nu}\Bigr)^{1 / \nu},
\]
where the infimum is taken over all possible couplings $(X, Y)$ such that the marginal distributions coincide, i.e., $X \deq X_0$, $Y \deq Y_0$. 

Similar to the Kolmogorov distance, Wasserstein distances can also characterize the convergence of distribution function: for $\nu\ge 1$ and a sequence of random variables $S_n$, $\lim_{n \to \infty} d_{\nu}(S_{n}, Y)=0 $ implies $S_n\dto Y$ \cite[Lemma~8.3]{Bickel1981}. More importantly, Wasserstein $1$-distance implies an upper bound on the Kolmogorov distance based on the proposition below. 
\begin{proposition}\label{prop:dist}
Suppose $X, Y$ are random variables and $Y$ is a stable random variable whose characteristic function is parametrized by \eqref{eq:GK} with a certain index $\alpha$. Then there exists a positive constant $B_\alpha$ such that
\[
\sup_{x\in \mathbb{R}}|\mathbb{P}(X\le x) - \mathbb{P}(Y\le x)| \le B_\alpha [d_1(X, Y)]^{1/2}.
\]
\end{proposition}
\begin{proof}[Proof of Proposition~\ref{prop:dist}]
By Kantorovich duality \cite{Rachev1984}, we have
\[
d_1(X, Y) =\sup_{h \in \mathfrak{L}}|\mathbb{E} h(Y)- \mathbb{E} h(X)|
\]
with $\mathfrak{L} \coloneqq \{h: \mathbb{R} \to \mathbb{R}:|h(y)-h(x)| \le |y-x|\}$.

Assume $h_0\in \mathfrak{L}$ and $h_0:\mathbb{R}\to[0, 1]$ such that $h_0(s) = 1$ for $s\le 0$, $h_0(s) = 0$ for $s\ge 1$, and $h_0(s)$ is a linear segment on $s\in(0, 1)$. Fix any $x\in\mathbb{R}$, and define $h(s) = \frac{1}{\psi}h_0(\psi(s - x))$ for some constant $\psi>0$. Then for any $y, z \in \mathbb{R}$, we have
\[
|h(y)-h(z)|=\frac{1}{\psi}\bigl|h_{0}\bigl(\psi(y-x)\bigr)-h_{0}\bigl(\psi(z-x)\bigr)\bigr| \le|y-z|,
\]
which implies $h\in\mathfrak{L}$. Then
\begin{align*}
\mathbb{P}(X \le x) & \le \mathbb{E}[h_{0}(\psi(X-x))]=\psi \mathbb{E}[h(X)] \\
& \le \psi \mathbb{E}[h(Y)]+\psi d_{1}(X, Y) \\
& \le \mathbb{P}\Bigl(Y \le x+\frac{1}{\psi}\Bigr) + \psi d_{1}(X, Y)
\end{align*}
Furthermore, by the boundedness of the stable densities, we can find a constant $B_\alpha$ such that $\mathbb{P}(Y \le x+\frac{1}{\psi})-\mathbb{P}(Y \le x) \le B_\alpha/\psi$,
which implies
\[
\mathbb{P}(X \le x)-\mathbb{P}(Y \le x) \le \psi d_{1}(X, Y)+B_\alpha \frac{1}{\psi}.
\]
Hence, we take $\psi=[d_1(X, Y)]^{-1/2}$ and have
\[
\mathbb{P}(X \le x)-\mathbb{P}(Y \le x) \le B_\alpha[d_1(X, Y)]^{\frac{1}{2}}.
\]
This gives one-half of the proposition. The other half follows similarly.
\end{proof}
The results above reveal the effectiveness of Wasserstein distance in terms of quantifying the rate of convergence of distribution functions, and thus many works have studied central limit theorems in Wasserstein distances. In particular, a relevant result in \cite{Johnson2005} is used in our proofs, which establishes a general rate of convergence (in Wasserstein $\nu$-distance) to stable distributions for sums of independent random variables. We attach the result below for completeness.
\begin{theorem}[Johnson--Samworth]\label{thm:JS}
Fix $\alpha \in(0,2)$, and let $X_1, X_2, \dots$ be independent random variables (where $\mathbb{E} X_{i}=0$ if $\alpha>1$), and $S_{n}=(X_{1}+\dots+X_{n})/ n^{1 / \alpha}$. If there exists a stable random variable $Y$ such that $\sup_{i} d_{\nu}(X_{i}, Y)<\infty$ for some $\nu \in(\alpha, 2]$, then $\lim_{n \to \infty} d_{\nu}(S_{n}, Y)=0$ and
\[
d_{\nu}(S_{n}, Y) \le \frac{2^{1 / \nu}}{n^{1 / \alpha}}\Bigl(\sum_{i=1}^{n} d_{\nu}^{\nu}(X_{i}, Y)\Bigr)^{1 / \nu}.
\]
\end{theorem}

\section{Main proofs}\label{sec:proof}
\subsection{Proofs in Section~\ref{sec:char}}\label{appx:char}
\begin{proof}[Proof of Theorem~\ref{thm:1char}]
Since the sufficiency has been derived in Corollary~\ref{cor:Cparam}, here we only show the necessity. Assume the distributional equation \eqref{eq:distri} holds for $p = 1$. We have
\[
 \frac{1}{X_1} \deq \frac{1}{c}\Bigl(\frac{1}{X_1} + \frac{1}{X_2}\Bigr).
\]
By iterating, we have
\[
 \frac{1}{X_1} \deq \frac{1}{c^2}\Bigl(\frac{1}{X_1}+\frac{1}{X_2}
 +\frac{1}{X_3}+\frac{1}{X_4}\Bigr)
 \deq \lim_{n\to\infty}\frac{1}{c^n} \sum_{i = 1}^{2^n}\frac{1}{X_i}
\]
where $X_i$, $i = 1, 2, \dots$ are i.i.d.\ random variables that have the same distribution as $X_1$. Hence, $1/X_1$ is the limiting distribution of $\frac{1}{n}\sum_{i = 1}^{n}X_i^{-1}$ given conditions \eqref{it:1char1}, \eqref{it:1char2}, and it belongs to stable distributions.

Next we will figure out the parameters $\alpha$ and $\beta$ in the canonical representation \eqref{eq:GK}. Denote by $G$ the distribution function of $1/X_i$. Conditions \eqref{it:1char1}, \eqref{it:1char2} imply that as $x\to\infty$,
\[
 \frac{G(-x)}{1 - G(x)}\to 1,
\]
and for every constant $k>0$,
\[
 \frac{1 - G(x) + G(-x)}{1 - G(kx) + G(-kx)}
 \to k.
\]
By Theorem~\ref{thm:GKclt}, we know $c_1 = c_2>0$, $\alpha = 1$. Then $\beta = (c_1 - c_2)/(c_1 + c_2) = 0$. The characteristic function $\varphi$ of $1/X_i$ has the expression
\[
 \log \varphi(t)=i \gamma t-\lambda|t|
\]
and satisfies
\[
 \varphi(t) = \varphi(t/c)^2.
\]
Solving the equation above gives $c = 2$. Finally, notice that the stable distribution with $\log \varphi(t) = i\gamma t - \lambda|t|$ is Cauchy distribution, i.e., $\Cau(\gamma, \lambda)$. Hence, $X_i$ still follows a Cauchy distribution after taking the reciprocal.
\end{proof}

\begin{proof}[Proof of Theorem~\ref{thm:2char}]
Since the sufficiency has been derived in Proposition~\ref{prop:2infdivsym}, again, we only show the necessity. By a similar argument as in the proof of Theorem~\ref{thm:1char}, we have $1/X_i^p\sim \Sta(r/p, -1, \lambda, 0)$.

An additional step is to show any distribution solving \eqref{eq:distri} is symmetric. By Assumption~\eqref{it:2char2}, we have
\[
 \lim _{n \to \infty} \mathbb{P}\bigl(\sgn(X_1X_2\cdots X_{n})= -1 \bigr)=\frac{1}{2}.
\]
By generalizing the result of \cite[Theorem~2.1]{Galambos2003}, we can show that $X_1$ is symmetric, i.e., $X_1 \deq -X_1$.
\end{proof}

\begin{proof}[Proof of Theorem~\ref{thm:infchar}]
Since the sufficiency has been derived in Proposition~\ref{prop:infparam}, again, we only show the necessity. Assume we have the equation
\[
 \frac{1}{X_1} \deq \frac{1}{c}\max\Bigl(\frac{1}{X_1}, \frac{1}{X_2}\Bigr)
\]
hold. Again, by iterating, we have
\[
 \frac{1}{X_1} \deq \frac{1}{c^2}\max\Bigl(\frac{1}{X_1}, \frac{1}{X_2}
 ,\frac{1}{X_3}, \frac{1}{X_4}\Bigr)\\
 \deq \lim_{n\to\infty}\frac{1}{c^n} \max_{i = 1, \dots, 2^n}\Bigl(\frac{1}{X_i}\Bigr).
\]
Denote by $G$ the distribution function of $1/X_i$. By Assumption~\eqref{it:infchar2}, we know $G(x) < 1$ for $x<\infty$. Assumption~\eqref{it:infchar1} in the theorem is equivalent to say for $k>0$,
\[
 \lim_{x\to\infty}\frac{1 - G(kx)}{1 - G(x)} = k^{-r}.
\]
By \cite[Theorem~2.1.1]{Galambos1978}, the asymptotic distribution for $\max_{i = 1, \dots, 2^n}(1 / X_i)$ under proper normalizing constants has the distribution function in the form of \eqref{eq:frechet}. Hence, $X_1$ should satisfy $X_1>0$, $X_1^r \sim \Exp(\lambda)$ for some $\lambda>0$, otherwise it would be degenerate. Expressing the distributional equation in  $G$ gives
\[
 G(x) = G(cx)^2.
\]
Since $G(x)$ has the form in \eqref{eq:frechet}, $c = 2^{1/r}$ solves the equation.
\end{proof}

\subsection{Proofs in Section~\ref{sec:clt}}\label{appx:clt}
Recall that $G$ is the distribution function of $1/X_i$ when $p = 1, \infty$; it is the distribution function of $1/X_i^p$ when $p$ is an even integer. 
\begin{proof}[Proof of Theorem~\ref{thm:1clt}]
For notational convenience, write
\[
 b_{n}=\inf \Bigl\{x: \mathbb{P}(|1 / X_1| >x) \le \frac{2}{n\pi} \le \mathbb{P}(|1 / X_1|\ge x)\Bigr\}.
\]
First consider the limiting distribution of $\frac{1}{b_n}\sum_{i = 1}^n X_i^{-1} - a_n$ with suitable $a_n$ and $b_n$. Observe that 
\[
 G(x) = \mathbb{P}\Bigl(\frac{1}{X_1}\le x\Bigr) =
 \begin{cases}
 F(0) + 1 - F(1/x^-) &\text{ if } x>0,\\
 F(0) - F(1/x^-) &\text{ if } x< 0,
 \end{cases}
\]
and $G(0) = \lim_{x \to 0^+}G(x) = F(0)$. As a result, \eqref{eq:1clt} implies \eqref{eq:clt1}. In addition, since $F(x) - F(-x)$ is regularly varying as $x\to 0$ with index $r = 1$, we have as $x\to\infty$
\begin{align*}
    \frac{1-G(x)+G(-x)}{1-G(k x)+G(-k x)} &=
    \frac{F(1/x) - F(-1/x^-)}{F(1/kx) - F(-1/kx^-)}\to k.
\end{align*}
By Theorem~\ref{thm:GKclt}, we know $\alpha = 1$ and $c_1 = c_2$, and therefore
\[
    \frac{1}{b_n}\sum_{i = 1}^n X_i^{-1} - a_n \dto \Cau(\gamma, \tau).
\] 
By choosing $a_n$, $b_n$ according to Theorem~\ref{thm:GKclt} and assigning $c_1 = c_2 = \frac{1}{\pi}$, we have $\gamma = 0$ and $\tau = 1$. 

Now, we analyze the limiting behavior of $a_n$ and $b_n$.
Given \eqref{eq:1clt3}, $\lim_{n\to\infty} \frac{n}{b_n} = \tau$. For $a_n$, notice that
\[
a_n = \frac{n}{b_n}\int_{|x|<b_n} x \, d G(x) =\frac{n}{b_n}\int_{|y|>1 / b_n}\frac{1}{y}\, dF(y).
\]
Based on \eqref{eq:1clt2} and \eqref{eq:1clt3}, we have $\lim_{n\to\infty} a_n = \tau\gamma$. 

Finally, we get back to establishing the asymptotic result for Shepp $1$-product. We define two random variables $C\sim\Cau\bigl(\gamma\tau^2 / (\gamma^2\tau^2 + 1), \tau / (\gamma^2\tau^2 + 1)\bigr)$, $C_0\sim\Cau(0, 1)$, respectively. For $x>0$, as $n\to\infty$, 
\begin{align*}
 \mathbb{P}\bigl(nS_1(X_1, \dots, X_n)\le x\bigr)
 &= \mathbb{P}\bigl(n S_1(X_1, \dots, X_n)< 0\bigr)
 + \mathbb{P}\bigl(0 < nS_1(X_1, \dots, X_n)\le x\bigr) \\
 &=\mathbb{P}\Bigl(\frac{1}{n}\sum_{i = 1}^n X_i^{-1} < 0\Bigr) + \mathbb{P}\Bigl(\frac{1}{n}\sum_{i = 1}^n X_i^{-1}\ge \frac{1}{x}\Bigr) \\
 &\to \mathbb{P}\Bigl(C_0< -\lim_{n\to\infty}a_n\Bigr)
 + \mathbb{P}\Bigl(C_0\ge \frac{1}{x}\lim_{n\to\infty}\frac{n}{b_n}-\lim_{n\to\infty}a_n\Bigr) \\
 &= \mathbb{P}(C_0< -\tau\gamma)
 + \mathbb{P}\Bigl(C_0\ge \frac{\tau}{x} - \tau \gamma\Bigr) \\
 &= \mathbb{P}(C \le x)
\end{align*}
where the weak convergence is due to Lemma~\ref{lemma:cdf}. The same result can be derived for $x<0$ analogously.
\end{proof}

\begin{proof}[Proof of Theorem~\ref{thm:2clt}]
We decompose our proof into two parts. First, we derive the asymptotic distribution of $1/S_p^p(X_1, \dots, X_n)$. Second, we show that the asymptotic distribution corresponding to $S_p(X_1, \dots, X_n)$ must be symmetric. These two observations uniquely determine the asymptotic distribution for Shepp $p$-product. 

Again, we consider the limiting distribution of $\frac{1}{b_n}\sum_{i = 1}^n X_i^{-p} - a_n$ with suitable $a_n$ and $b_n$. Observe that 
\[
 G(x) = \mathbb{P}\Bigl(\frac{1}{X_1^p}\le x\Bigr) =
 \begin{cases}
 F(-1/x^{1/p}) + 1 - F({1/x^{1/p}}^-) &\text{ if } x>0,\\
 0 &\text{ if } x\le 0.
 \end{cases}
\]
As a result, \eqref{eq:clt1} naturally holds with $c_1 = 0$. In addition, since $F(x) - F(-x)$ is regularly varying as $x\to 0$ with index $r = p\alpha$, we have as $x\to\infty$
\begin{equation}\label{eq:2da2}
    \frac{1-G(x)+G(-x)}{1-G(k x)+G(-k x)} \to k^{\alpha}.
\end{equation}
By Theorem~\ref{thm:GKclt}, we know the limiting distribution is a stable distribution with $\alpha\in(0, 1)$ and $c_1 = 0$. By choosing $a_n$, $b_n$ according to Theorem~\ref{thm:GKclt} and assigning $c_2 = 1 / \bigl(\Gamma(1-\alpha)\cos(\pi\alpha/2)\bigr)$, we have $\beta = -1$, $\lambda = 1$, and $\gamma = 0$.

Note we can choose $a_n = 0$ due to Theorem~\ref{thm:GKclt}. Given \eqref{eq:2clt}, we have $\lim_{n\to\infty} \frac{n^{1/\alpha}}{b_n} = \sigma$. We define two random variables $S_0\sim\Sta(\alpha, -1, 1, 0)$, $S_\sigma\sim \Sta(\alpha, -1, \sigma^{-\alpha}, 0)$, respectively. As $n\to\infty$, we have for $x>0$,
\begin{align*}
 \mathbb{P}\Bigl(\frac{1}{n^{1/\alpha}}\frac{1}{S_p^p(X_1, \dots, X_n)}\le x\Bigr)
 &= \mathbb{P}\Bigl(\frac{1}{n^{1/\alpha}}\sum_{i = 1}^{n}\frac{1}{X_i^p}\le x\Bigr)\\
 &\to \mathbb{P}\Bigl(S_0\le x \lim_{n\to\infty}\frac{n^{1/\alpha}}{b_n}\Bigr)\\
 &= \mathbb{P}(S_\sigma\le x )
\end{align*}
where the weak convergence is due to Lemma~\ref{lemma:cdf}. Same result can be derived for $x<0$ analogously.

It remains to show the symmetry of the limiting distribution, which can be established by \cite[Theorem~2.1]{Galambos2003}. The regular variation condition in the theorem implies that $F(x)$ is continuous at zero, and therefore $\mathbb{P}(X_i = 0)=0$.
Since $F(0)\in(0, 1)$, we have $\mathbb{P}(X_i< 0)\in(0, 1)$, $\mathbb{P}(X_i > 0)\in (0,1)$, and thus
\[
 \lim_{n \to\infty} \mathbb{P}\bigl(\sgn(X_1X_2\cdots X_{n})=-1\bigr)=\frac{1}{2}.
\]
Given the asymptotic distribution of $1/S_p^p(X_1, \dots, X_n)$, we can generalize the result of \cite[Theorem~2.1]{Galambos2003} to show that $S_p(X_1, \dots, X_n)$ enjoys a symmetric limiting distribution. It is then uniquely determined as in the theorem. 
\end{proof}

\begin{proof}[Proof of Theorem~\ref{thm:infclt}]
We consider the limiting distribution of $\max (1/X_1,\dots,1/X_n)/b_n$ with suitable $b_n$. According to \cite[Theorem~2.1.1]{Galambos1978}, necessary and sufficient conditions for $\max (1/X_1,\dots,1/X_n)/b_n$ to converge to an extreme value distribution with the distribution function \eqref{eq:frechet} are that
\begin{equation}\label{eq:infda}
 \sup\{x:G(x)<1\} = \infty
\end{equation}
and that for all $x>0$,
\begin{equation}\label{eq:infda2}
 \lim_{t\to \infty}\frac{1 - G(tx)}{1 - G(t)} = x^{-r}.
\end{equation}
Under conditions \eqref{eq:infda} and \eqref{eq:infda2}, we can choose $b_{n}=\inf \{x: 1-G(x)\le1 / n\}$ and have
\[
 \frac{\max (1/X_1,\dots,1/X_n)}{b_n} \dto
 U_0^{-1}
\]
where $U_0$ is a positive random variable such that $U_0^r \sim \Exp(1)$.

Notice that for $x>0$, $1 - G(x) = \mathbb{P}(1/X_1 >x) = \mathbb{P}(0<X_1<1/x) = F(1/x^-) - F(0)$. Hence, \eqref{eq:infda} is satisfied under the assumption that $F(x)>F(0)$ for $x>0$. Condition \eqref{eq:infda2} is also satisfied by regular variation condition in the theorem. We express $b_n$ by
\[
 b_n = \inf \{x: F(1 / x)-F(0) \le 1 / n\}.
\]
The condition \eqref{eq:infclt} implies that $\lim_{n\to \infty} b_n/n^{1/r} = \lambda$.

Now, we get back to establishing the asymptotic result for Shepp $\infty$-product. For $x>0$, as $n\to\infty$,
\begin{align*}
 \mathbb{P}\bigl(n^{1/r} S_\infty(X_1,\dots,X_n) \le x\bigr)
 &= \mathbb{P}\Bigl(\frac{\max (1/X_1,\dots,1/X_n)}{n^{1/r}} \le 0 \Bigr)
 + \mathbb{P}\Bigl(\frac{\max (1/X_1,\dots,1/X_n)}{n^{1/r}} \ge \frac{1}{x}\Bigr) \\
 &= \mathbb{P}\Bigl(\frac{\max (1/X_1,\dots,1/X_n)}{b_n} \le 0\Bigr)
 + \mathbb{P}\Bigl(\frac{\max (1/X_1,\dots,1/X_n)}{b_n} \ge \frac{n^{1/r}}{b_n x}\Bigr) \\
 &\to \mathbb{P}(U_0^{-1} \le 0)
 + \mathbb{P}\Bigl( U_0^{-1} \ge \frac{1}{x}\lim_{n\to \infty}\frac{n^{1/r}}{b_n}\Bigr)
 = \mathbb{P}\bigl( U \le x\bigr)
\end{align*}
where the weak convergence is due to Lemma~\ref{lemma:cdf}. For $x\le 0$, we observe that as $n\to\infty$,
\begin{align*}
 \mathbb{P}\bigl(n^{1/r} S_\infty(X_1,\dots,X_n) \le x\bigr)
 &\le \mathbb{P}(\max (1/X_1,\dots,1/X_n) \le 0) \\
 &\le \mathbb{P}(X_1<0)^n \to 0.
\end{align*}
Since $F(x)>F(0)$ for $x>0$, we have $\mathbb{P}(X_1<0) < 1$ and the limit above holds.
\end{proof}

\begin{proof}[Proof of Theorem~\ref{thm:1rate}]
We first develop an equivalent expression for our target Kolmogorov distance. Then we can simply apply the results in \cite{Hall1981}.

Let $C_\tau$ be a random variable such that $C_\tau \sim \Cau(0, \tau)$ and $C$ be the random variable defined in the theorem. By the transformation property of the Cauchy distribution,
\[
\frac{1}{\gamma+C_{1/\tau}}
\sim
\Cau\Bigl(
\frac{\gamma \tau^2}{\gamma^2\tau^2+1},
\frac{\tau}{\gamma^2\tau^2+1}
\Bigr).
\]
Thus it holds that
\[
C\deq\frac{1}{\gamma+C_{1/\tau}}.
\]
For $x<0$, we have
\begin{align*}
 |F_{n}(x)-\mathbb{P}( C\le x)|
 &=\biggl|\mathbb{P}\biggl(1/x\le \frac{\sum_{i = 1}^n 1/X_i}{n} < 0 \biggr)-\mathbb{P}\bigl( 1/x - \gamma \le C_{1/\tau} < - \gamma \bigr)\biggr| \\
 &\le \biggl|\mathbb{P}\biggl(\frac{\sum_{i = 1}^n 1/X_i}{n} - \gamma \le -\gamma \biggr)-\mathbb{P}\bigl( C_{1/\tau} \le - \gamma \bigr)\biggr| \\
 &+ \biggl|\mathbb{P}\biggl(\frac{\sum_{i = 1}^n 1/X_i}{n} - \gamma \le \frac{1}{x}-\gamma \biggr)-\mathbb{P}\Bigl( C_{1/\tau} \le \frac{1}{x} - \gamma \Bigr)\biggr| \\
 &\le 2\sup_{x\in\mathbb{R}}|G_n(x) - G_0(x)|
\end{align*}
where we denote by $G_n$ the distribution function of $\frac{1}{n}\sum_{i = 1}^n 1/X_i - \gamma$, and $G_0$ the distribution function of $C_{1/\tau}$. For $x>0$, the same bound can be obtained analogously. Hence,
\[
 \sup_{x\in\mathbb{R}}|F_{n}(x)-\mathbb{P}( C\le x)| \le 2\sup_{x\in\mathbb{R}}|G_n(x) - G_0(x)|.
\]

Consider the sum $(\sum_{i = 1}^n1/X_i - n\gamma)/ n$ and its limiting distribution $C_{1/\tau}$. We verify the conditions in \cite[Corollary~2]{Hall1981} to derive the final result. First, by assumption, $S(x)$ and $D(x)$ are ultimately monotone, therefore are ultimately of the one sign. Second, $D(x)$ is absolutely integrable as $|D(x)| = O(x^{-2}\log^2 x)$. This verifies the condition in \cite[Corollary~2]{Hall1981} and completes the proof.
\end{proof}

\begin{proof}[Proof of Theorem~\ref{thm:2rate}]
Since $F$ is symmetric and continuous (thus, $F$ has no point mass), for $z > 0$, we have
\begin{align*}
 \mathbb{P}\bigl(n^{1/(p\alpha)} S_p(X_1, \dots, X_n) \le z\bigr)&=\frac{1}{2}+\frac{1}{2} \mathbb{P}\biggl(\frac{n^{1/(p\alpha)}}{(\sum_{i=1}^{n} 1 / X_i ^p )^{1/p}} \le z\biggr) \\
 &=\frac{1}{2}+\frac{1}{2} \mathbb{P}\biggl(\frac{\sum_{i=1}^{n} 1 / X_i ^p }{n^{1/\alpha} } \ge \frac{1}{z^p }\biggr).
\end{align*}
Define $G$ as the distribution function of $1/X_i^p$ and a random variable $S_0 \sim \Sta(\alpha, -1, \sigma^{-\alpha}, 0)$. As in Theorem~\ref{thm:1rate}, we have
\begin{align}
 \sup_{x\in\mathbb{R}}|F_{n}(x)-F_0(x)| &=
 \frac{1}{2}\sup_{x>0}\biggl|\mathbb{P}\biggl(\frac{\sum_{i=1}^{n} 1 / X_i ^p }{n^{1/\alpha} } \le \frac{1}{x^p }\biggr)-\mathbb{P}\Bigl( S_0 \le \frac{1}{x^p }\Bigr)\biggr| \nonumber\\
 &= \frac{1}{2} \sup_{x\in\mathbb{R}}|G_{n}(x)-G_0(x)|, \label{eq:2dist}
\end{align}
where $G_n$ is the distribution function of the sum $(\sum_{i=1}^{n} 1 / X_i ^p ) / n^{1/\alpha} $, and $G_0$ is the distribution function of $S_0$, i.e., the stable law to our interest.
Under the notation in \cite{Hall1981}, we have
\begin{align*}
 S(x)&=1-G(x)+G(-x)-(c_2 +c_1 ) \sigma^{-\alpha} x^{-\alpha}, \\
 D(x)&=1-G(x)-G(-x)-(c_2 -c_1) \sigma^{-\alpha} x^{-\alpha},
\end{align*}
 where $c_1$, $c_2$ are parameters in the canonical representation \eqref{eq:GK}. After some rearrangements, we have $S(x) = D(x) = a_1(x)x^{-\alpha}$.

In the following, we apply \cite[Corollary~3]{Hall1981} to derive \eqref{eq:2rate2}. Since $S(x)$ is continuous and $x^\varepsilon |S(x)|$ is ultimately nonincreasing, we know that $S(x)$ is ultimately monotone, and thus the conditions in \cite[Corollary~3]{Hall1981} are satisfied. Therefore,
\[
 \sup_{x\in\mathbb{R}}|G_{n}(x)-G_0(x)|=O(n^{1-s/\alpha})\quad \text{ as } \quad n\to \infty
\]
if and only if
\begin{equation}\label{eq:2iff}
 |G(x)-G_0(x)|+|G(-x)-G_0(-x)|=O(x^{-s}) \quad \text{ as }\quad x\to \infty.
\end{equation}
Recall that by \eqref{eq:2dist}, two Kolmogorov distances are equivalent. Besides, \eqref{eq:2iff} can be reduced to
\begin{equation}\label{eq:2iff2}
\bigl|F(x^{-1/p}) - F(-x^{-1/p})
-\mathbb{P}(S_0>x)\bigr|=O(x^{-s}), \quad s\in (\alpha, \min(1, 2\alpha))
\end{equation}
by the nonnegative support of $G$, $G_0$. For $S_0$ which follows a stable distribution $\Sta(\alpha, -1, \sigma^{-\alpha}, 0)$, it is shown in \cite[Equation~2.2]{Mijnheer1986} that as $x\to\infty$,
\[
\mathbb{P}(S_0 > x) = c_2 \sigma^{-\alpha} x^{-\alpha} + O(x^{-2\alpha}).
\]
Since $s<2\alpha$, Equation \eqref{eq:2iff2} is satisfied when
\[
|S(x)| = O(x^{-s}).
\]
Finally, a direct application of \cite[Corollary~3]{Hall1981} gives the claimed bound.
\end{proof}

\begin{proof}[Proof of Theorem~\ref{thm:2rate2} under continuous cdfs]
Denote by $G_i$, $G_0$ the distribution function of $1/X_i^2$, $1/Z_\sigma^2$ for $i = 1, 2, \dots$.
It suffices to give an explicit bound for
$\sup_{x\in\mathbb{R}}|F_{n}(x)-\mathbb{P}(Z_\sigma\le x)|$, for which we have
\[
    \sup_{x\in\mathbb{R}}|F_{n}(x)-\mathbb{P}
    (Z_\sigma\le x)| = \frac{1}{2} \sup_{x\in\mathbb{R}} |G_n(x)-G_0(x)|
\]
in Theorem~\ref{thm:2rate}.
By Proposition~\ref{prop:dist}, we further deduce that
\[
    \sup _{x\in\mathbb{R}}\bigl|G_{n}(x)-G_0(x)\bigr|
    \le \frac{\sigma^2}{2}\Bigl[d_1\Bigl(\frac{\sum 1/X_i^2}{n^2}, \frac{1}{Z_\sigma^2}\Bigr)\Bigr]^{1/2}
\]
since the constant $B_\alpha \le \sigma^2/2$.

In the following, we derive a concrete bound for the Wasserstein-1 distance term above using Theorem \ref{thm:JS}. First, for each $i$, we apply the identity for the Wasserstein-1 distance \cite{Villani2021topics} to obtain
\begin{equation}\label{eq:w1}
d_1\Bigl(\frac{1}{X_i^2}, \frac{1}{Z_\sigma^2}\Bigr) = \int_{\mathbb{R}} |G_i(x) - G_0(x)| \, dx. 
\end{equation}

For $x>0$, by the continuity of $F_i$ and the assumption in the theorem, we have
\[
    1-G_i(x)
    = \mathbb{P}(|X_i|<x^{-1/2}) = \bigl(c_2+a_i(x)\bigr)x^{-1/2}.
\]
On the other hand,
\[
    1-G_0(x)
    =\mathbb{P}(|Z_\sigma| < x^{-1/2})
    =2 \int_0^{\frac{1}{\sigma\sqrt{x}}} \frac{1}{\sqrt{2\pi}} e^{-t^2 / 2} \, dt.
\]
For $u\ge 0$, one can verify that 
\[
    0\le \sqrt{\frac{2}{\pi}}u - 2 \int_0^{u} \frac{1}{\sqrt{2\pi}} e^{-t^2 / 2} \, dt
    =\sqrt{\frac{2}{\pi}} \int_0^u(1-e^{-t^2/2})\,dt
    \le \frac{u^3}{3\sqrt{2\pi}}.
\]
Consequently, by setting $u = 1 / (\sigma\sqrt{x})$, we obtain
\[
    |1-G_0(x)-c_2 x^{-1/2}|
    \le C_0x^{-3/2},
    \quad C_0\coloneqq\frac{1}{3\sigma^3\sqrt{2\pi}}.
\]
It follows that, for every $x\ge 1$ and every $i$,
\begin{align*}
    |G_i(x)-G_0(x)| \le C_i x^{-1/2-\gamma_i}+C_0 x^{-3/2} \le Mx^{-1/2-\gamma_*}+C_0x^{-3/2}.
\end{align*}
In the last inequality, we adopt the notation that $M \coloneqq \sup_{i \in \mathbb{N}_+} C_i < \infty$ and $\gamma_* \coloneqq \inf_{i \in \mathbb{N}_+} \gamma_i > 1/2$. Then, by combining the analysis above and Equation \eqref{eq:w1}, we obtain
\begin{align*}
    d_1(\frac{1}{X_i^2},\frac{1}{Z_\sigma^2}) &\le \int_0^1 |G_i(x) - G_0(x)| \, dx + \int_1^\infty |G_i(x) - G_0(x)| \, dx \\
    &\le 1 + M \int_1^\infty x^{-1/2-\gamma_*}\,dx
        + C_0 \int_1^\infty x^{-3/2}\,dx \\
    &= 1+\frac{M}{\gamma_*-1/2} + 2C_0.
\end{align*}
In the last equation, we use the fact that $\gamma_* > 1/2$. Also, note that the upper bound above holds for any $i$.

Lastly, applying the Wasserstein-1 bound to Theorem \ref{thm:JS}, we obtain
\begin{equation}\label{eq:2wass-rate-alternative}
    d_1\Bigl(\frac{\sum 1/X_i^2}{n^2}, \frac{1}{Z_\sigma^2}\Bigr)
    \le \frac{2}{n^2}\sum_{i=1}^n d_1(\frac{1}{X_i^2},\frac{1}{Z_\sigma^2}) = O\Bigl(\frac{1}{n}\Bigr).
\end{equation}
This completes the proof.
\end{proof}

\begin{proof}[Proof of Proposition~\ref{prop:infdna}]
We apply \cite[Theorem~2.1.1]{Galambos1978} to confirm the convergence.
For a sequence $\max(X_1,\dots, X_n)/b_n$, necessary and sufficient conditions for it to converge to an extreme value distribution are that
\[
 \sup\{x:G(x)<1\} = \infty
\]
and that for all $x>0$,
\[
 \lim_{t\to \infty}\frac{1 - G(tx)}{1 - G(t)} = x^{-r}.
\]
Since $\lambda^r + a_1(x)\to \lambda^r$ as $x \to \infty$ and $a>0$, the first condition above is satisfied. Plugging in the expression of $G$ verifies the second condition. Hence, for i.i.d.\ random variables $X_i$ that follow $G$, we have
\[
 \frac{\max_i(X_1, \dots, X_n)}{b_n} \dto \frac{1}{\lambda U}
\]
and $b_n$ is given by
\[
 b_n = \inf \{x: 1 - G(x) \le 1 / n\}.
\]
It is easy to check that $\lim b_n/n^{1/r} = \lambda$. Hence, by Lemma~\ref{lemma:cdf}, our proof is complete.
\end{proof}

\begin{proof}[Proof of Theorem~\ref{thm:infrate}]
We prove the case $\lambda = 1$, as the general settings can be derived up to a change of scale. For $x\le0$, we have $F_{r}(x) = 0$ and
\[
 \mathbb{P}\Bigl(\frac{\max_i(X_1, \dots, X_n)}{n^{1/r}}\le x\Bigr)\le\mathbb{P}\Bigl(\frac{\max_i(X_1, \dots, X_n)}{n^{1/r}}\le 0\Bigr) = G^n(0).
\]
Since $G(0)<1$, we can find a constant $C_1$, free of $n$, such that $G^n(0)\le C_1/n^s$ for all $n$. Next, we consider the case $x>0$, for which we rely on the mean value theorem.

We use $b_n = n^{1/r}$ throughout the proof and define $H(x) = e^{-x}$, $y_1 = -n\log G(b_n x)$, and $y_2 = x^{-r}$. We have
\begin{align*}
 H(y_1) &= G^n(b_n x) = \mathbb{P}\Bigl(\frac{\max_i(X_1, \dots, X_n)}{n^{1/r}}\le x\Bigr), \\
 H(y_2) &= F_{r}(x), \\
 H(y_1) - H(y_2) &= (y_1 - y_2) H^\prime(\xi),
\end{align*}
where $\xi$ is between $y_1$ and $y_2$. Hence, $|\mathbb{P}\bigl(\max_i(X_1, \dots, X_n)/ n^{1/r}\le x\bigr) - F_r(x)| = |y_1 - y_2| |H^\prime(\xi)|$.

Let $t = \bigl(1 + a_1(b_n x)\bigr)(b_n x)^{-r}$. By Taylor expansion,
\begin{equation}\label{eq:infrate}
 y_1 = -n\log (1 - t)
 = -n \Bigl(-t - \frac{t^2}{2(1 - \theta t)^2}\Bigr)
 = nt + \frac{nt^2}{2(1 - \theta t)^2}
\end{equation}
for some $\theta\in[0, 1]$. Let $\delta_0$ be a positive constant such that 
\begin{enumerate}[\upshape (i)]
    \item\label{it:infrate1} $a_1(x) \ge 0$ for $x\ge \delta_0$,
    \item\label{it:infrate2} $\bigl(1 + a_1(x)\bigr)/{x^{r}} \le \frac{1}{2}$ for $x\ge \delta_0$. 
    \end{enumerate}
Under the assumptions in the theorem, $\delta_0$ is well-defined. 

First, we consider $b_n x>\delta_{0}$, in which case $(1 - \theta t)^2 \ge 1/4$ by \eqref{it:infrate2}. With that being said, $(1 - \theta t)^2$ is bounded away from zero. To analyze $|y_1 - y_2||H^\prime(\xi)|$, note that
\begin{equation}\label{eq:infrate2}
 t = \frac{1 + a_1(b_n x)}{(b_n x)^{r}}
 = \frac{1}{nx^r} + \frac{a_1(b_n x)}{nx^r}
\end{equation}
and
\[
 \frac{|a_1(b_n x)|}{n(x)^r}
 \le \frac{C_0}{n^{1 + \eta/r}(x)^{r + \eta}}
\]
since $|a_1(x)|\le C_0/|x|^\eta$. We have
\begin{align*}
 |y_1 - y_2| & \le \biggl|nt + \frac{nt^2}{2(1 - \theta t)^2} - x^{-r}\biggr| \\
 & \le \frac{C_0}{n^{\eta/r}(x)^{r + \eta}} + 2n \biggl(\frac{1}{nx^r} + \frac{a_1(b_n x)}{n(x)^r}\biggr)^2 \\
 & \le \frac{C_0}{n^{\eta/r}(x)^{r + \eta}} + 4n \biggl(\frac{1}{n^2x^{2r}} + \frac{a_1^2(b_n x)}{n^2(x)^{2r}}\biggr) \\
 & \le \frac{C_0}{n^{\eta/r}(x)^{r + \eta}} + \frac{4}{n x^{2r}} + \frac{4C_0^2}{n^{1 + 2\eta/r}(x)^{2r + 2\eta}}.
\end{align*}

For $|H^\prime(\xi)| = e^{-\xi}$, we have
\[
 |H^\prime(\xi)| \le \max(e^{-y_1}, e^{-y_2})
 \le e^{-y_2} \max( e^{y_2 - y_1}, 1).
\]
From \eqref{eq:infrate} and \eqref{eq:infrate2}, we have
\[
 y_2 - y_1 = x^{-r} - nt - \frac{nt^2}{2(1 - \theta t)^2}
 = -\frac{a_1(b_n x)}{(x)^r} - \frac{nt^2}{2(1 - \theta t)^2}.
\]
Noticing that $y_2 - y_1\le 0$ by \eqref{it:infrate1}, we get $|H^\prime(\xi)| \le e^{-x^{-r}}$ and
\begin{align}
 |y_1 - y_2||H^\prime(\xi)|
 &\le \biggl( \frac{C_0}{n^{\eta/r}(x)^{r + \eta}} + \frac{4}{n x^{2r}} + \frac{4C_0^2}{n^{1 + 2\eta/r}(x)^{2r + 2\eta}}\biggr) e^{-x^{-r}} \nonumber\\
 &\le \underbrace{%
 \biggl(\frac{C_0}{(x)^{r + \eta}} + \frac{4}{ x^{2r}} + \frac{4C_0^2}{(x)^{2r + 2\eta}} \biggr) e^{-x^{-r}}
 }_{c(x)} \cdot \frac{1}{n^s}. 
 \label{eq:infineq}
\end{align}
Now, we need to upper bound $c(x)$ using some constant $C_2$. Since $\lim_{x\to 0^+}c(x) = \lim_{x\to \infty} c(x) = 0$ and $c(x)$ is continuous on $(0, \infty)$, there exists a finite constant $C_2$ such that $c(x) \le C_2$ for all $x$, $b_n x > \delta_0$.

Next, we consider $x>0$ while $b_n x\le \delta_0$, in which case \eqref{eq:infineq} does not hold. We will show that both $G^n(b_n x)$ and $F_r(x)$ decay at an exponential rate, and hence can be bounded by a polynomial $C_3/n^s$. First, we have $G^n(b_n x)\le G^n(\delta_{0})$ with $G(\delta_{0})<1$. Second, $F_r(x)\le e^{- n \delta_{0}^{-r}}$ due to $x\le\delta_{0}/ n^{1/r}$. Therefore, there exists $C_3$ such that $|\mathbb{P}\bigl(\max_i(X_1, \dots, X_n)/ n^{1/r}\le x\bigr) - F_r(x)| \le C_3/n^s$ for all $x$, $b_n x\le \delta_0$.

To sum up, by choosing $C = \max(C_1, C_2, C_3)$, we have
\[
\sup_{x\in\mathbb{R}}\bigl|\mathbb{P}\Bigl(\frac{\max_i(X_1, \dots, X_n)}{n^{1/r}}\le x\Bigr) - F_r(x) \bigr| \le \frac{C}{n^{s}}.
\]
\end{proof}

\subsection{Proofs in Section~\ref{sec:test}}\label{appx:test}

\begin{proof}[Proof of Theorem~\ref{thm:null}]
We first derive the asymptotic distribution for $\mathit{IU}_n^E$. Note that $\mathit{IU}_n^E$ is asymptotically equivalent to the $U$-statistic with a centered kernel of degree three,
\[
\Psi_E(x, y, z)=\frac{1}{3}(\mathbb{I}_{\{k_E(x, y)<z\}}+\mathbb{I}_{\{k_E(x, z)<y\}}+\mathbb{I}_{\{k_E(y, z)<x\}})-\frac{1}{2}.
\]
Then we calculate the projection of the kernel $\Psi_E$. Given $x$, let
\begin{align*}
\psi_E(x) &\coloneqq  \mathbb{E}\bigl(\Psi_E(X, Y, Z)\bigm| X = x\bigr) \\
&=\frac{1}{3}\bigl[\mathbb{P}\bigl(k_E(x, Y)<Z\bigr)+\mathbb{P}\bigl(k_E(x, Z)<Y\bigr)+\mathbb{P}\bigl(k_E(Y, Z)<x\bigr)\bigr]-\frac{1}{2}.
\end{align*}

Without loss of generality, we assume
$X_i \oversim{iid} \Exp(1)$ under $H_0$, and $F$ is the distribution function of $\Exp(1)$.
By the characterization result, $\mathbb{P}\bigl(k_E(Y, Z)<x\bigr) = F(x)$. Since $X, Y, Z$ are i.i.d.\ random variables, $\mathbb{P}\bigl(k_E(x, Y)<Z\bigr) = \mathbb{P}\bigl(k_E(x, Z)<Y\bigr)$ and for $x>0$
\begin{align*}
\mathbb{P}\bigl(k_E(x, Y)<Z\bigr) &= \mathbb{P}\biggl(\max\Bigl(\frac{1}{2x}, \frac{1}{2Y}\Bigr)>\frac{1}{Z}\biggr) \\
&= \mathbb{P}\Bigl(\frac{1}{2Y}>\frac{1}{Z}\Bigr) + \mathbb{P}\Bigl(\frac{1}{2x}>\frac{1}{Z}\Bigr) - \mathbb{P}\Bigl(\frac{1}{2x}>\frac{1}{Z}, \frac{1}{2Y}>\frac{1}{Z}\Bigr) \\
&= \frac{1}{3} + \frac{2}{3}e^{-3x}.
\end{align*}
Hence, we have
\[
\psi_E(x) = \frac{1}{3}F(x) + \frac{2}{3}\mathbb{P}\bigl(k_E(x, Y)<Z\bigr) - \frac{1}{2} = \frac{1}{18} - \frac{1}{3}e^{-x} + \frac{4}{9}e^{-3x}
\]
and $\var{\psi_E(X)} = \frac{11}{3780}$. By \cite[Theorem~7.1]{Hoeffding1948},
\[
\sqrt{\frac{420n}{11}}\mathit{IU}_n^E \dto \Nor(0, 1).
\]

Next, we consider Kolmogorov-type statistic. For a fixed $t>0$, the difference $U_n^E(t) - F_n(t)$ is same as a $U$-statistic with kernel
\[
\Xi_E(x, y ; t)=\mathbb{I}_{\{k_E(x, y)<t\}}-\frac{1}{2}(\mathbb{I}_{\{x<t\}} + \mathbb{I}_{\{y<t\}}).
\]
Now, the projection of the kernel is
\[
\xi_E(x ; t)\coloneqq \mathbb{E}\bigl(\Xi_E(X, Y ; t) \bigm| X=x\bigr)=\mathbb{P}\bigl(k_E(x, Y)<t\bigr)-\frac{1}{2}\bigl(\mathbb{I}_{\{x<t\}} + F(t)\bigr).
\]
For $\mathbb{P}\bigl(k_E(x, Y)<t\bigr)$ with $x>0$, $t>0$, we have the expression as follows
\[
\mathbb{P}\bigl(k_E(x, Y)<t\bigr) =
\begin{cases}
1 & \text{ if } t>2x, \\
F(t/2) & \text{ if } t\le 2x.
\end{cases}
\]
Therefore,
\[
\xi_E(x ; t) =
\begin{cases}
1/2 -F(t)/2 & \text{ if } t>2x, \\
F(t/2) - F(t)/2 - \mathbb{I}_{\{x<t\}}/2& \text{ if } t\le 2x.
\end{cases}
\]
Now, we calculate the variance of $\xi_E(X; t)$ and get
\[
\var\bigl(\xi_E(X; t)\bigr) =
\frac{1}{4}e^{-t}(1 - e^{-t}).
\]
The variance is positive for $t>0$, indicating a nondegenerate kernel.
Using the methods of the work \cite{Silverman1983}, we can show that the U-statistical empirical process $\eta_n(t) = \sqrt{n}(U_n^E(t) - F_n(t))$ converges weakly as $n \to \infty$ to some centered Gaussian process $\eta(t)$ with covariance
\[
\mathbb{E} \eta(s)\eta(t) = 4 \mathbb{E} (\xi_E(X; s)\xi_E(X; t)), \qquad s<t.
\]
\end{proof}

\subsection{Proofs in Section~\ref{sec:app}}\label{appx:app}
\begin{proof}[Proof of Theorem~\ref{thm:align1}]
    We denote all the words by $w_1, w_2, \dots, w_{d^k}$, and define the covariance matrix $\Sigma_0$ for words $w_1, w_2, \dots, w_{d^k-1}$ according to \cite[Equation~5]{Reinert2009}. From \cite{Reinert2009}, we have
    \[
        (X_{w_1}, \dots, X_{w_{d^k - 1}}, Y_{w_1}, \dots, Y_{w_{d^k - 1}})\dto(Z_1, Z_2)
    \]
    under the mild conditions that $m\ge n$, $k<\frac{\bar{n}}{2}$. The word counts $X_{w_{d^k}}, Y_{w_{d^k}}$ are omitted because we always have 
    $X_{w_{d^k}} = -\sum_{j = 1}^{d^k - 1} X_{w_j}$, $Y_{w_{d^k}} = -\sum_{j = 1}^{d^k - 1} Y_{w_j}$. Then, by the continuous mapping theorem, we have $D_p^S\dto D_p^{\lim}$ for $p = 1, \infty$ and even numbers.
    \end{proof}

\end{document}